\documentclass[amssymb, notitlepage,amsmath, pra, 10pt, aps]{revtex4-1}
\usepackage{amsmath,amssymb,amsfonts} 
\usepackage{graphicx}
\usepackage{bm}
\usepackage{psfrag}
\usepackage[caption=false]{subfig}
\usepackage{color}
\usepackage{xcolor}
\usepackage{hyperref}
\usepackage{mathtools}

\DeclareMathOperator*{\argmin}{arg\,min}

\newtheorem{theorem}{Theorem}
\newtheorem{lemma}{Lemma}
\newtheorem{definition}{Definition}
\newtheorem{proof}{Proof}

\allowdisplaybreaks[1]

\begin{document}

\title{Distributed Computing for Huge-Scale Linear Programming}

\author{Luoyi Tao}
\email[]{luoyitao@smail.iitm.ac.in, taoluoyi@gmail.com}
\affiliation{Geophysical Flows Lab\\
                and\\
             Department of Aerospace Engineering\\
             Indian Institute of Technology Madras\\
             Chennai 600 036, India}

\date{\today}

\def\s{\!}
\def\ss{\!\!}
\def\sss{\!\!\!}
\def\l{\left}
\def\r{\right}
\def\Prox{\text{Prox}}
\def\IndicatorFunction{\delta}
\def\ConeK{{\cal K}}
\def\SubjectTo{\text{subject to}}
\def\cA{{\cal A}}
\def\cB{{\cal B}}
\def\cC{{\cal C}}
\def\cD{{\cal D}}
\def\cE{{\cal E}}
\def\cF{{\cal F}}
\def\cG{{\cal G}}
\def\cH{{\cal H}}
\def\cI{{\cal i}}
\def\cK{{\cal K}}
\def\cL{{\cal L}}
\def\cP{{\cal P}}
\def\cQ{{\cal Q}}
\def\cR{{\cal R}}
\def\dist{\textbf{dist}}
\def\dom{\textbf{dom}}
\def\Indicator{I}
\def\maximize{\text{maximize}}
\def\minimize{\text{minimize}}
\def\overX{\bar{X}}
\def\overY{\bar{Y}}
\def\Projection{\Pi}
\def\prox{\textbf{prox}}
\def\subjectto{\text{subject to}}
\def\Lik{\Lambda_i^k}
\def\Ljk{\Lambda_j^k}
\def\NCV{N_{\text{CV}}}
\def\LDmu{\mu}
\def\LDnu{\nu}
\def\LDbeta{\beta}
\def\pcE{E}
\def\pcF{F}
\def\pcG{G}
\def\pcH{H}
\def\dvE{\mu}
\def\dvF{\lambda}
\def\dvG{\mu}
\def\dvH{\nu}

\def\dvZ{\mu} 
\def\dvZp{\prescript{p}{}{\!\dvZ}}
\def\dvZn{\prescript{n}{}{\!\dvZ}}
\def\svZ{Y}
\def\svZp{\prescript{p}{}{\!\svZ}}
\def\svZn{\prescript{n}{}{\!\svZ}}
\def\svZpcoef{\prescript{p}{}{\!\gamma}}
\def\svZncoef{\prescript{n}{}{\!\gamma}}

\def\dvH{\mu} 
\def\dvHp{\prescript{p\!H}{}{\!\dvH}}
\def\dvHn{\prescript{n\!H}{}{\!\dvH}}
\def\svH{Y}
\def\svHp{\prescript{p\!H}{}{\!\svH}}
\def\svHn{\prescript{n\!H}{}{\!\svH}}
\def\svHpcoef{\prescript{p\!H}{}{\!\gamma}}
\def\svHncoef{\prescript{n\!H}{}{\!\gamma}}

\def\dvZ{\mu} 
\def\svZ{Y}
\def\dvZP{\prescript{p\!}{}{\!\dvZ}}
\def\dvZN{\prescript{n\!}{}{\!\dvZ}}
\def\svZP{\prescript{p\!}{}{\!\svZ}}
\def\svZN{\prescript{n\!}{}{\!\svZ}}
\def\svZPcoef{\prescript{p\!}{}{\!\gamma}}
\def\svZNcoef{\prescript{n\!}{}{\!\gamma}}

\def\dvH{\mu} 
\def\dvHP{\prescript{p\!H}{}{\!\dvH}}
\def\dvHN{\prescript{n\!H}{}{\!\dvH}}
\def\svH{Y}
\def\svHP{\prescript{p\!H}{}{\!\svH}}
\def\svHN{\prescript{n\!H}{}{\!\svH}}
\def\svHPcoef{\prescript{p\!H}{}{\!\gamma}}
\def\svHNcoef{\prescript{n\!H}{}{\!\gamma}}

\begin{abstract}
This study develops an algorithm for distributed computing of linear programming problems of huge-scales. Global consensus with single common variable, multiblocks, and augmented Lagrangian are adopted. The consensus is used to partition the constraints of equality and inequality into multi-consensus blocks, and the subblocks of each consensus block are employed to partition the primal variables into $M$ sets of disjoint subvectors. The global consensus constraints of equality and other constraints are replaced equivalently by the extended constraints of equality involving slack variables, since the slack variables help the feasibility and  initialization of the algorithm. The block-coordinate Gauss-Seidel method, the proximal point method, and ADMM are used to update the primal variables, descent models used to update the dual. Convergence of the algorithm to optimal solutions is argued and the rate of convergence, $O(1/k^{1/2})$ is estimated, under feasibility of the algorithm and boundedness of the dual sequences supposed. Analysis is presented on how to ensure the feasibility and boundedness through initial and control parameter values and a dual descent model with built-in bound for the original constraints of inequality. Further exploration of dual descent models with built-in bound is needed.
\begin{description}
\item[Mathematics Subject Classification]
\verb+90C05+
\verb+90C06+
\verb+90C25+
\verb+90C30+
\verb+68W15+
\verb+68W40+
\verb+76F02+
\verb+76F05+
\verb+76F10+
\end{description}
\end{abstract}

\maketitle

\section{Introduction}

The problem studied here arose naturally from the exploration of hydrodynamic turbulence modeling
within the framework of optimal correlation functions \cite{Tao2020}.
The formulation involves a huge number of the multi-point correlations of 
turbulent velocity and pressure fluctuations
as control variables
and a huge number of constraints of equality and inequality, affine and
quadratic, generated by the Navier-Stokes equations,
  the Cauchy-Schwarz inequality, the non-negativity of variance of products,
and $\overline{(a-c)^2}\geq 0$ where $a$ and $c$ denote fluctuations and the overline
the ensemble averaging.
The basic idea and technical details in the context of homogeneous shear turbulence
can be found in \cite{Tao2020}, which presents some preliminary numerical simulation results
and is yet to be updated (with enlarged number of control variables,
errors corrected,
the total turbulent energy as an objective that is maximized,
and some new numerical tests).
The difficulty in testing the framework adequately lies in that the algorithms and software available
like conventional ADMM and scs 
\cite{Boydetal2010, ODonoghueetal2016, ODonoghueetalPackage2016c, ParikhBoyd2013} 
cannot deal with the computational size of the problem and
are restricted to convex constraints,
 it is essential to develop distributed computing algorithms
for huge-scale optimizations.

Considering the challenge in the development of such algorithms,
we restrict our study to linear programming (LP) in this work.
There are several reasons for this restriction:
One is that the homogeneous shear turbulence problem may be approximated
as linear programming.
The other is that the study may provide insights on how to extend 
the methodology to certain nonlinear programming
and CFD-optimization coupled problems relevant to turbulence modeling.
The third is that huge-scale linear programming independently deserves
special attention.

The study exploits conventional techniques of optimization 
\cite{BC2011, Beck2017, Boydetal2010}
along with newer developments
and organizes them in a rather natural manner, 
both underlying justifications and mathematical analysis are presented.
(i) The global consensus with single common variable is applied to
partition the constraints of equality and inequality into multiple consensus blocks;
the subblocks of each consensus block further partition the primal variables
into $M$ sets of disjoint subvectors
in order to make computation of the primal variables feasible. 
These consensus constraints of equality and other constraints
are converted equivalently into the extended constraints of equality involving slack variables;
the slack variables provide the flexibility to initialize the algorithm and
 make it feasible.
(ii) The augmented Lagrangian is adopted as the basis to solve for the primal
variables iteratively;
the block-coordinate Gauss-Seidel method, the proximal point method, and ADMM are
 used to update the primal variables.
 The penalty terms in the augmented Lagrangian help bound the dual sequences.
 (iii) The descent models for updates of the dual variables are introduced 
 as a simple and natural tool to satisfy the constraints by the limits
 of iterations. Here, the dual update rules are viewed as
 algorithmic modeling analogous to constitutive modeling in continuum mechanics,
 which may offer us certain flexibility in developing algorithms for optimization.
 (iv) Under the algorithm supposedly feasible and the dual sequences supposedly bounded,
 convergence of the algorithm to optimal solution is argued 
 and the rate of convergence, $O(1/k^{1/2})$ is estimated.
 (v) The feasibility and boundedness needs to be ensured through 
 appropriate initialization of the primal and dual sequences
and choice of the control parameter values.
To help resolve the issue of boundedness,
  we extend the dual model for the original constraints of inequality
  to the dual model with built-in bound.
 
The paper is organized as follows.
Section~\ref{sec:Formulation} formulates the optimization problem.
Section~\ref{sec:UpdatePrimal} presents methodology for updates of the primal variables.
The descent updates for the dual variables and the analysis of convergence
of the algorithm are discussed in Sec.~\ref{sec:UpdateDual}.
The issues of feasibility of the algorithm,
boundedness of the dual sequences, initialization, control parameter values,
and precondition of the constraints
are discussed in some detail in Sec.~\ref{sec:Initialization}.
Section~\ref{sec:Summary} summarizes the main results and lists some issues
that need to be explored further.

\section{Formulation of optimization problem}
\label{sec:Formulation}

The primal problem is linear programming,
\begin{align}
\label{LPGSPrimalProblemOriginal}
&
\minimize\ \ \  f(Z)
\notag\\[-6.5pt]&
\\[-6.5pt]&
\subjectto\ \ Z\in\cC.
\notag
\end{align}
Here, $f(Z)$ is linear; 
\begin{align}
 \cC=\{Z\in \mathbb{R}^n: 
 \pcG(Z)\leq 0,\ \pcH(Z)=0,\ l\leq Z\leq u\},
\label{LPGSPrimalProblemConstraintSeto}
\end{align}
$\pcG$ and $\pcH$ are vector functions whose components are affine;
$l$ and $u$ are the lower and upper bounds for $Z$, respectively.
Without loss of generality, we take $u>0$ and $l=-u$,
which can be obtained through a translational shift of the domain,
$f$ can be treated equivalently as linear and $\pcG$ and $\pcH$ as affine in the shifted domain.
Further, to have adequate initialization of the algorithm to be developed,
the constraints of equality are converted equivalently into the constraints of inequality,
$\{\pcH(Z)\leq 0$, $-\pcH(Z)\leq 0\}$ and included as part of $\pcG(Z)\leq 0$,
\begin{align}
 \cC=\{Z\in [-u, u]: \pcG(Z)\leq 0\}.
\label{LPGSPrimalProblemConstraintSet}
\end{align}

Considering the huge-scale computational size,
both the value of $n$ and the number of constraints contained 
in \eqref{LPGSPrimalProblemConstraintSet} being great,
we first partition \eqref{LPGSPrimalProblemOriginal} constraint-wise
via global consensus with single common variable $Z$ 
\cite{Boydetal2010, ParikhBoyd2013},
\begin{align}
\label{LPGSPrimalProblemGlobalConsensusRing}
&
 \minimize\ \ \, \sum_{i=1}^N \big(f(X_i) +\Indicator_{\cC_i}(X_i)\big)
\notag\\[-13pt]&
\\[-2pt]&
\subjectto\ \ X_i-Z=0, \ \  i=1,\ldots,N.
\notag
\end{align}
Here, $\cC$ is partitioned into $N$ blocks, $\cC=\cup_{i=1}^{N}\cC_i$
with $\cC_i$ given by 
\begin{align}
 \cC_i=\{X_i\in [-u, u]:
    \pcG_i(X_i)\leq 0\},
\label{LPGSConstraintsithCBo}
\end{align}
and $\Indicator_{\cC_i}$ is the indicator function of $\cC_i$. 
Also, to have appropriate initialization of the algorithm, 
$\{X_i-Z=0\}$ is replaced equivalently by $\{Z-X_{i}\leq 0$, $X_i-Z\leq 0\}$
and $\cC_i$ of \eqref{LPGSConstraintsithCBo} is extended to
\begin{align}
 \cC_i=\{X_i\in [-u, u]:
        Z-X_{i}\leq 0,\ X_i-Z\leq 0,\  Z\in [-u,u],\
    \pcG_i(X_i)\leq 0\}.
\label{LPGSConstraintsithCB}
\end{align}

To help initialize the algorithm, we introduce the slack variables
so as to convert all the constraints of inequality in $\cC_i$ 
of \eqref{LPGSConstraintsithCB}
into the extended constraints of equality,
 following the conventional practice
(and viewing the slack variables as part of the primal variables below),
\begin{align}
\cC_i=\big\{
        X_i\in [-u, u]:\,&
        Z-X_{i}+\svZp_{i}=0,\ X_{i}-Z+\svZn_{i}=0,\
        \svZp_{i}\in[0, u_{\svZp_{i}}],\ 
        \svZn_{i}\in[0, u_{\svZn_{i}}],\ 
        Z\in [-u,u],\
        \notag\\&
        \pcG_i(X_i)+Y_i= 0,\ 
        Y_{i}\in[0, u_{Y_i}]
    \big\}.
\label{LPGSConstraintsithCBSlacked}
\end{align}
The upper bounds
for the slack variables 
$\{u_{\svZp_{i}}, u_{\svZn_{i}}, u_{Y_{i}}\}$
are set as follows.
(a) $u_{\dvZp_{i}}=u_{\dvZn_{i}}=2u+ \epsilon_{Z}$ for some $\epsilon_{Z}\geq 0$.
(b) $u_{Y_{i}}=\sum_{l=1}^{M}\sum_{j}\vert \pcG_{i,l,j}\vert\, u_{l,j}+\pcG_{i,0}
        +\epsilon_{\pcG_{i}}$ for some $\epsilon_{\pcG_{i}}\geq 0$. 
Here, 
$\pcG_{i}(Z)$ $=\sum_{l=1}^{M}\pcG_{i,l} Z_{l}+\pcG_{i,0}
        =\sum_{l=1}^{M}\sum_{j}\pcG_{i,l,j}Z_{l,j}+\pcG_{i,0}$ is used.
For convenience, we call this $\cC_i$ the $i$-th consensus block (CB)
or the $i$-CB to indicate the relevant operations involved. 
                
The above treatment involves the conversion of the global consensus constraints of equality
in \eqref{LPGSPrimalProblemGlobalConsensusRing} to
the equivalent inequality constraints in \eqref{LPGSConstraintsithCB}
and further to the extended constraints of equality involving slack variables
in \eqref{LPGSConstraintsithCBSlacked}.
Though increasing computational size of the optimization problem,
it has the advantage of providing required initialization
 to make the algorithm feasible because the slack variables can be initialized rather flexibly, 
as to become clear later.

To help solve the objective function of \eqref{LPGSPrimalProblemGlobalConsensusRing}
subject to \eqref{LPGSConstraintsithCBSlacked},
 we employ the augmented Lagrangian $L$,
\vspace{-3mm}
\begin{align}
 L(X,Z,\svZp, \svZn, Y, \dvZp, \dvZn, \dvG, \rho)
=\sum_{i=1}^N L_i(X_i,Z, \svZp_{i}, \svZn_{i}, Y_{i}, \dvZp_{i}, \dvZn_{i}, \dvG_{i}, \rho_i).
\label{LPGSLagranginDualConsensus}
\end{align}
Here,
\begin{align*}
&
X=(X_1,\ldots,X_N),\ 
\svZp=(\svZp_1,\ldots, \svZp_N),\
\svZn=(\svZn_1,\ldots, \svZn_N),\
Y=(Y_1,\ldots, Y_N),\
\dvZp=(\dvZp_1,\ldots, \dvZp_N),
\notag\\&
\dvZn=(\dvZn_1,\ldots, \dvZn_N),\
\dvG=(\dvG_1,\ldots, \dvG_N),\
\rho=(\rho_1,\ldots, \rho_N),
\end{align*}
(the transpose symbol is ignored to avoid cumbersome notation),
and the augmented Lagrangian function for the $i$-CB, $L_{i}$ takes the form of
\begin{align}
L_i(X_i,Z, \svZp_{i},\svZn_{i}, Y_{i}, \dvZp_{i}, \dvZn_{i}, \dvG_{i}, \rho_i)
=\,&
 f(X_{i})
+\langle  \dvZp_{i}, Z_{}-X_{i}+\svZp_{i}\rangle
+\langle  \dvZn_{i}, X_{i}-Z_{}+\svZn_{i}\rangle
+\langle  \dvG_{i}, \pcG_{i}(X^{}_{i})+Y_{i}\rangle
\notag\\&
+\frac{\rho_{i}}{2}\Big(
             \Vert Z_{}-X_{i}+\svZp_{i} \Vert^2
            +\Vert X_{i}-Z_{}+\svZn_{i} \Vert^2
            +\Vert \pcG_{i}(X^{}_{i})+Y_{i} \Vert^2
            \Big),
\label{LPGSLagranginFunctionithBC}
\end{align}
where $\rho_i$ is the positive penalty parameter to augment the Lagrangian function,
$\{\dvZp_{i}$, $\dvZn_{i}$, $\dvG_{i}\}$
are the dual variables associated with
$\{Z_{}-X_{i}+\svZp_{i}=0$, $X_{i}-Z_{}+\svZn_{i}=0$, $\pcG_{i}(X^{}_{i})+Y_{i}=0\}$, respectively.

The dual problem of \eqref{LPGSPrimalProblemGlobalConsensusRing} subject to 
\eqref{LPGSConstraintsithCBSlacked} is
\begin{align}
 \min_{\tilde{\cC}} \, \sup_{\dvZp,\,\dvZn,\,\dvG}
                        L(X,Z,\svZp, \svZn, Y, \dvZp, \dvZn,\dvG,\rho),
\label{LPGSDualProblemGlobalConsensusRing}
\end{align}
with 
$ \tilde{\cC}:= \{
        (X_{i},Z,\svZp_{i},\svZn_{i},Y)
        \in[-u, u]\times[-u, u]\times[0,u_{\svZp_{i}}]\times[0, u_{\svZn_{i}}]
        \times[0, u_{Y_{i}}],\ i=1,\ldots, N
\}$.
Suppose that the problem~\eqref{LPGSDualProblemGlobalConsensusRing} has a saddle-point 
$(X^{\ast},Z^{\ast}$, $\dvZp^{\ast}$,
$\dvZn^{\ast}$, $\dvG^{\ast})$
without the penalty and slack variables,
  $\{\rho=0$, $\svZp=\svZn=0$, $Y=0\}$.
 We have the saddle point theorem,
\begin{align}
 L(X^{\ast},Z^{\ast},0,0,0,\dvZp,\dvZn,\dvG, 0) 
 \leq L(X^{\ast},Z^{\ast},0,0,0,\dvZp^{\ast},\dvZn^{\ast},\dvG^{\ast},0)
 \leq L(X,Z,0,0,0,\dvZp^{\ast},\dvZn^{\ast},\dvG^{\ast},0),
\label{LPGSSaddlePointTheorem}
\end{align}
and the associated Karush-Kuhn-Tucker (KKT) conditions, $\forall i=1,\ldots, N,$
\vspace{-3mm}
\begin{align}
&
\dvZp^{\ast}_{i}-\dvZn^{\ast}_{i}
 \in \partial_{X_{i}}
        \Big[    f(X^{\ast}_{i})
                +\langle  \dvG^{\ast}_{i}, \pcG_{i}(X^{\ast}_{i})\rangle
                +\Indicator_{[-u,u]}(X^{\ast}_{i})
        \Big],\ \
-\frac{1}{N}\sum_{i=1}^{N}(\dvZp^{\ast}_{i}-\dvZn^{\ast}_{i})
\in \partial_{Z}\Indicator_{[-u,u]}(Z^{\ast}),\ 
\notag\\[-3pt]&
 X^{\ast}_{i}=Z^{\ast}_{},\ 
 \pcG_{i}(X^{\ast}_{i})\leq 0,\ 
 \dvG^{\ast}_{i}\geq 0,\ 
 \langle \dvG^{\ast}_{i}, \pcG_{i}(X^{\ast}_{i}) \rangle=0.
\label{LPGSKKTConditions}
\end{align}

To make it computationally feasible,  $X_i$ and $Z$ are partitioned 
into $M$ disjoint subvectors, respectively,
\begin{align}
X_{i}=(X_{i,1},\ldots, X_{i,M}),\
Z=(Z_{1},\ldots, Z_{M}),\ \forall i=1,\ldots,N,
\label{LPGSXiZiPartitioned}
\end{align}
where the dimensions and  component orders of subvectors
$\{X_{i,l}, Z_{l}\}$ are independent of $i$;
$X_{i,l}\in\mathbb{R}^{m_l}$, $\sum_{l=1}^{M}m_l=n$.
It is supposed that $m_{l}$ for all $l$ have similar values.
To help solve  $X_{i,l}$  through an iterative procedure,
 we introduce 
\begin{align}
&
X^{j,+}_{i,l}
:=(X^{j}_{i,1},\ldots, X^{j}_{i,l-1}, X_{i,l}, X^{k}_{i,l+1}, \ldots, X^{k}_{i,M}),\
X^{j,k}_{i,l}
:=(X^{j}_{i,1},\ldots, X^{j}_{i,l-1}, X^{k}_{i,l}, X^{k}_{i,l+1}, \ldots, X^{k}_{i,M}),
\notag\\&
X^{j,k+1}_{i,l}
:=(X^{j}_{i,1},\ldots, X^{j}_{i,l-1},
        X^{k+1}_{i,l}, X^{k}_{i,l+1}, \ldots, X^{k}_{i,M}),\ 
X^{k}_{i}:=(X^{k}_{i,1},\ldots, \ldots, X^{k}_{i,M}),\
X^{k}:=\{X^{k}_{i}: \forall i\},
\label{LPGSConstraintsithCBPartitioned}
\end{align}
where the superscripts $j, k$ denote the $j/ k$-th iteration
and the subscript $(i,l)$ denotes  the subblock
or  the process to be updated.
(We could use a more clear but longer notation like $X^{j,+,k}_{i,l}$.)
In the definitions of \eqref{LPGSConstraintsithCBPartitioned},
$X_i$ can be substituted by $Z$, $\svZp_{i}$, and $\svZn_{i}$,
respectively. 
Accordingly, $\cC_i$ of \eqref{LPGSConstraintsithCBSlacked} is partitioned 
by a block coordinate Guess-Seidel method,
\vspace{-1mm}
\begin{align}
 \cC_i=\cup_{l=1}^{M}\cC_{i,l},\
 \cC_{i,l}:=\big\{\,&
        Z-X^{k+1,+}_{i,l}+\svZp_{i}=0,\ X^{k+1,+}_{i,l}-Z_{}+\svZn_{i}=0,\
        \pcG_i(X^{k+1,+}_{i,l})+Y_i= 0, 
        \notag\\&
        \svZp_{i}\in[0, u_{\svZp_{i}}],\ 
        \svZn_{i}\in[0, u_{\svZn_{i}}],\ 
        Y_{i}\in[0, u_{Y_i}],\ 
        Z\in [-u,u],\
        X_{i,l}\in [-u_{l}, u_{l}]
    \big\}.
\label{LPGSConstraintsithCBlthSubblock}
\end{align}
These equality constraints 
may be violated by $X^{k+1,+}_{i,l}$ at finite $k$, owing to the approximate 
nature of iterative algorithm; for this very reason
 the augmented Lagrangian method is applied since it allows for violation 
of the constraints in intermediate steps of iteration;
the constraints are satisfied in the limit of $k\rightarrow\infty$
and the constraints~\eqref{LPGSConstraintsithCBlthSubblock} need to be viewed as such.
The above partition is suitable to distributed computing
and is computational in nature on the basis of iterations.
Next, $L_i$ of \eqref{LPGSLagranginFunctionithBC} 
is partitioned functionally according to
\vspace{-2mm}
\begin{align}
&
 L_i
=\sum_{l=1}^{M}L_{i,l}(X^{k+1,+}_{i,l},Z_{l}, \svZp_{i,l}, \svZn_{i,l}, Y_{i}, 
                        \dvZp_{i,l}, \dvZn_{i,l}, \dvG_{i}, \rho_i),
\notag\\&
L_{i,l}(X^{k+1,+}_{i,l},Z_{l}, \svZp_{i,l}, \svZn_{i,l}, Y_{i},
                        \dvZp_{i,l}, \dvZn_{i,l}, \dvG_{i}, \rho_i)
\notag\\&
=f(X^{k+1,+}_{i,l})
    +\langle  \dvZp_{i,l}, Z_{l}-X_{i,l}+\svZp_{i,l}\rangle
    +\langle  \dvZn_{i,l}, X_{i,l}-Z_{l}+\svZn_{i,l}\rangle
    +\langle  \dvG_{i}, \pcG_{i}(X^{k+1,+}_{i,l})+Y_{i}\rangle
\notag\\&\hskip5mm
+\frac{\rho_{i}}{2}\Big(
             \Vert Z_{l}-X_{i,l}+\svZp_{i,l} \Vert^2
            +\Vert X_{i,l}-Z_{l}+\svZn_{i,l} \Vert^2
            +\Vert \pcG_{i}(X^{k+1,+}_{i,l})+Y_{i} \Vert^2
            \Big).
\label{LPGSLagranginFunctionithBClthSubblock}
\end{align}
The roles of the penalty quadratic terms are to become clear
in Sec.~\ref{sec:Initialization}.

\section{Update Rules for Primal Variables}
\label{sec:UpdatePrimal}
We update the primal variables $X_{i,l}$ 
 by applying the block coordinate Gauss-Seidel method,
the proximal point method, and ADMM,
on the basis of \eqref{LPGSLagranginFunctionithBClthSubblock}.
At iteration $k$,
fix $i\in\{1,\ldots,N\}$ and $l\in \{1,\ldots,M\}$. The primal variable
 $X_{i,l}$  of the $(i,l)$-subblock is updated through
\vspace{-3mm}
\begin{align} 
X^{k+1}_{i,l}=
\argmin_{-u_{l}\leq X_{i,l}\leq u_{l}}
    \Big\{\,&
         L_{i,l}(X^{k+1,+}_{i,l},Z^{k}_{l}, \svZp^{k}_{i,l}, \svZn^{k}_{i,l}, Y^{k}_{i},
                    \dvZp^{k}_{i,l}, \dvZn^{k}_{i,l},
                        \dvG^{k}_{i}, \rho_i)
        +\frac{\sigma^{k}_{i}}{2}\Vert X_{i,l}-X^{k}_{i,l} \Vert^2
    \Big\},
\label{LPGSProximalXilA}
\end{align}
where $\sigma^{k}_{i}$ is a positive proximal control parameter for the proximal point algorithm
common to all the subblocks of the $i$-CB and its value is to be estimated 
in Sec.~\ref{sec:Initialization}.

Based on \eqref{LPGSLagranginFunctionithBClthSubblock},
$Z_{l}$ is updated through
\vspace{-3mm}
\begin{align}
Z^{k+1}_{l}=
\argmin_{-u_{l}\leq Z_{l}\leq u_{l}} 
        \sum_{i=1}^{N}\Big(\,&
             \langle  \dvZp^{k}_{i,l}, Z_{l}-X^{k+1}_{i,l}+\svZp^{k}_{i,l}\rangle
            +\langle  \dvZn^{k}_{i,l}, X^{k+1}_{i,l}-Z_{l}+\svZn^{k}_{i,l}\rangle
            \notag\\[-3pt]&
            +\frac{\rho_{i}}{2}\big(
                                     \Vert Z_{l}-X^{k+1}_{i,l}+\svZp^{k}_{i,l} \Vert^2
                                    +\Vert X^{k+1}_{i,l}-Z_{l}+\svZn^{k}_{i,l} \Vert^2
                                    \big)
            +\frac{\tau^{k}}{2}\Vert Z_{l}-Z^{k}_{l} \Vert^2\Big),
\label{LPGSProximalZl}
\end{align}
where $\tau^{k}$ is a positive proximal control parameter whose value is to be estimated
in Sec.~\ref{sec:Initialization}. 
Considering that this summation operation may be too big to be implemented in a single process,
 the following simple rudimentary procedure may be designed to realize it:
($\big(\sum_{j=1}^{0} (\tau^{k}+2\rho_{j})\big) Z_{0,l}$ is set to zero.)
\vspace{-2mm}
\begin{align}
&
\text{In the $(i,l)$-process, $i=1,\ldots,N-1$},
\notag\\&
\hskip10mm \text{$\Big(\sum_{j=1}^{i} (\tau^{k}+2\rho_{j})\Big) Z_{i,l}:=
\Big(\sum_{j=1}^{i-1} (\tau^{k}+2\rho_{j})\Big) Z_{i-1,l}
        +2\rho_{i}X^{k+1}_{i,l}
             +\rho_{i}(\svZn^{k}_{i,l}-\svZp^{k}_{i,l})+\dvZn^{k}_{i,l}-\dvZp^{k}_{i,l}
                        +\tau^{k}Z^{k}_{l}$,}
             \notag\\&\hskip50mm 
             \text{and pass to the $(i+1,l)$-process.}
\notag\\&
\text{In the $(N,l)$-process,}
\notag\\&\hskip10mm
\text{$\Big(\sum_{j=1}^{N} (\tau^{k}+2\rho_{j})\Big) Z_{N,l}:=
\Big(\sum_{j=1}^{N-1} (\tau^{k}+2\rho_{j})\Big) Z_{N-1,l}
        +2\rho_{N}X^{k+1}_{N,l}
             +\rho_{N}(\svZn^{k}_{N,l}-\svZp^{k}_{N,l})+\dvZn^{k}_{N,l}-\dvZp^{k}_{N,l}
                        +\tau^{k}Z^{k}_{l}$,} 
 \notag\\& \hskip10mm
\text{$Z^{k+1}_{l}:=P_{[-u_{l},u_{l}]}(Z_{N,l})$ and pass to the $(i,l)$-processes for all
                    $i=1,\ldots,N-1$.}
\label{LPGSProximalZkAlgorithm}
\end{align}
Here, the symbols $Z_{i,l}$ are used as intermediate variables 
to help presentation, with $i$ indicating the computation done in the $i$-CB,
$P_{[-u_{l},u_{l}]}$ is the projection into $[-u_{l},u_{l}]$.

Since all the slack variables hold for the whole $i$-CB,
 $X^{k+1}_{i}$ and $Z^{k+1}$ are used to update them.
Fix $i\in\{1,\ldots,N\}$.
\vspace{-3mm}
\begin{align} 
\svZp^{k+1}_{i}=
\argmin_{0\leq \svZp_{i}\leq u_{\svZp_{i}}}\s
\Big\{
 \langle \dvZp^{k}_{i}, Z^{k+1}_{}-X^{k+1}_{i}+\svZp_{i}\rangle
+\frac{\rho_{i}}{2}\Vert Z^{k+1}_{}-X^{k+1}_{i}+\svZp_{i} \Vert^2
+\frac{\svZpcoef^{k}_{i}}{2}\Vert \svZp_{i}-\svZp^{k}_{i} \Vert^2
    \Big\},
\label{LPGSProximalsvZpilA}
\end{align}
\vspace{-3mm}
\begin{align} 
\svZn^{k+1}_{i}=
\argmin_{0\leq \svZn_{i}\leq u_{\svZn_{i}}}\s
\Big\{
 \langle \dvZn^{k}_{i}, X^{k+1}_{i}-Z^{k+1}_{}+\svZn_{i}\rangle
+\frac{\rho_{i}}{2}\Vert X^{k+1}_{i}-Z^{k+1}_{}+\svZn_{i} \Vert^2
+\frac{\svZncoef^{k}_{i}}{2}\Vert \svZn_{i}-\svZn^{k}_{i} \Vert^2
    \Big\},
\label{LPGSProximalsvZnilA}
\end{align}
and
\begin{align} 
Y^{k+1}_{i}=
\argmin_{0\leq Y_{i}\leq u_{Y_{i}}}\s
\Big\{
 \langle \dvG^{k}_{i}, \pcG_{i}(X^{k+1}_{i})+Y_{i}\rangle
+\frac{\rho_{i}}{2} \Vert \pcG_{i}(X^{k+1}_{i})+Y_{i} \Vert^2
+\frac{\gamma^{k}_{i}}{2}\Vert Y_{i}-Y^{k}_{i} \Vert^2
    \Big\}.
\label{LPGSProximalYilA}
\end{align}
Here, $\svZpcoef^{k}_{i}$, $\svZncoef^{k}_{i}$, and $\gamma^{k}_{i}$
are positive proximal control parameters whose numerical values are to be estimated
in Sec.~\ref{sec:Initialization}.

\section{Dual Updates and Convergence Analysis}
\label{sec:UpdateDual}

To analyze the feasibility and convergence of the algorithm composed of the primal updates
\eqref{LPGSProximalXilA}, \eqref{LPGSProximalZl}, \eqref{LPGSProximalsvZpilA},
\eqref{LPGSProximalsvZnilA}, \eqref{LPGSProximalYilA},
and the dual updates~\eqref{LPGSDualUpdateRules} to be introduced,
we apply the first-order characterization of convex functions
to the functions involved in the primal updates.

First, the update \eqref{LPGSProximalXilA} gives
\begin{align} 
&
-\rho_{i}\Big(
             X^{k+1}_{i,l}-Z^{k}_{l}-\svZp^{k}_{i,l}
            +X^{k+1}_{i,l}-Z^{k}_{l}+\svZn^{k}_{i,l}
            +(\pcG_{i,l})^T\s(\pcG_{i}(X^{k+1,k+1}_{i,l})+Y^{k}_{i})
            \Big)
-\sigma^{k}_{i}(X^{k+1}_{i,l}-X^{k}_{i,l})
\notag\\
\in\,&
\partial_{X_{i,l}}
    \Big\{
         f(X^{k+1,+}_{i,l})
        +\big\langle  \dvZn^{k}_{i,l}-\dvZp^{k}_{i,l}, X_{i,l}\big\rangle
        +\big\langle  \dvG^{k}_{i}, \pcG_{i}(X^{k+1,+}_{i,l})\big\rangle
        +\Indicator_{[-u_{l}, u_{l}]}(X_{i,l})
    \Big\}_{X^{k+1}_{i,l}},
\label{LPGSProximalXilAFOC}
\end{align}
where $\Indicator_{[-u_{l}, u_{l}]}(X_{i,l})$ is the indicator function.
The first-order characteristic of the convex function involved on the right-hand side of
\eqref{LPGSProximalXilAFOC} under $X_{i,l}=X^{k}_{i,l}$ is
\begin{align} 
&
 \Indicator_{[-u_{l}, u_{l}]}(X^{k}_{i,l})
+f(X^{k+1,k}_{i,l})
-\big\langle \dvZp^{k}_{i,l}, X^{k}_{i,l} \big\rangle
+\big\langle \dvZn^{k}_{i,l}, X^{k}_{i,l} \big\rangle
+\big\langle \dvG^{k}_{i}, \pcG_{i}(X^{k+1,k}_{i,l}) \big\rangle
\notag\\&
+\frac{\rho_{i}}{2}\Big(\Vert Z^{k}_{l}-X^{k}_{i,l}+\svZp^{k}_{i,l} \Vert^2
                +\Vert X^{k}_{i,l}-Z^{k}_{l}+\svZn^{k}_{i,l} \Vert^2
                +\Vert \pcG_{i}(X^{k+1,k}_{i,l})+Y^{k}_{i} \Vert^2\Big)
\notag\\
\geq\,&
 \Indicator_{[-u_{l}, u_{l}]}(X^{k+1}_{i,l})
+f(X^{k+1,k+1}_{i,l})
-\big\langle \dvZp^{k}_{i,l}, X^{k+1}_{i,l} \big\rangle
+\big\langle \dvZn^{k}_{i,l}, X^{k+1}_{i,l} \big\rangle
+\big\langle \dvG^{k}_{i}, \pcG_{i}(X^{k+1,k+1}_{i,l}) \big\rangle
\notag\\&
+\frac{\rho_{i}}{2}\Big(\Vert Z^{k}_{l}-X^{k+1}_{i,l}+\svZp^{k}_{i,l} \Vert^2
                +\Vert X^{k+1}_{i,l}-Z^{k}_{l}+\svZn^{k}_{i,l} \Vert^2
                +\Vert \pcG_{i}(X^{k+1,k+1}_{i,l})+Y^{k}_{i} \Vert^2\Big)
\notag\\&
+\frac{\rho_{i}}{2}\Vert \pcG_{i}(X^{k+1,k+1}_{i,l})-\pcG_{i}(X^{k+1,k}_{i,l}) \Vert^2
+\big(\sigma^{k}_{i}+\rho_{i}\big)\Vert X^{k+1}_{i,l}-X^{k}_{i,l} \Vert^2.
\label{LPGSProximalXilAFOCSpecific}
\end{align}
Next, application of $\sum_{l=1}^{M}$ to \eqref{LPGSProximalXilAFOCSpecific},
use of $X^{k+1,k+1}_{i,l}=X^{k+1,k}_{i,l+1}$, 
and then operation of $\sum_{i=1}^{N}$ result in
\begin{align} 
&
 \sum_{i=1}^{N}\Big[ 
 \Indicator_{[-u_{}, u_{}]}(X^{k}_{i})
+f(X^{k}_{i})
-\big\langle \dvZp^{k}_{i}, X^{k}_{i} \big\rangle
+\big\langle \dvZn^{k}_{i}, X^{k}_{i} \big\rangle
+\big\langle \dvG^{k}_{i}, \pcG_{i}(X^{k}_{i}) \big\rangle
\notag\\[-3pt]&\hskip8mm
+\frac{\rho_{i}}{2}\Big(\Vert Z^{k}_{}-X^{k}_{i}+\svZp^{k}_{i} \Vert^2
                +\Vert X^{k}_{i}-Z^{k}_{}+\svZn^{k}_{i} \Vert^2
                +\Vert \pcG_{i}(X^{k}_{i})+Y^{k}_{i} \Vert^2\Big)
\Big]
\notag\\[-3pt]
\geq\,&
 \sum_{i=1}^{N}\Big[ 
 \Indicator_{[-u_{}, u_{}]}(X^{k+1}_{i})
+f(X^{k+1}_{i})
-\big\langle \dvZp^{k}_{i}, X^{k+1}_{i} \big\rangle
+\big\langle \dvZn^{k}_{i}, X^{k+1}_{i} \big\rangle
+\big\langle \dvG^{k}_{i}, \pcG_{i}(X^{k+1}_{i}) \big\rangle
\notag\\[-3pt]&\hskip8mm
+\frac{\rho_{i}}{2}\Big(\Vert Z^{k}_{}-X^{k+1}_{i}+\svZp^{k}_{i} \Vert^2
                +\Vert X^{k+1}_{i}-Z^{k}_{}+\svZn^{k}_{i} \Vert^2
                +\Vert \pcG_{i}(X^{k+1}_{i})+Y^{k}_{i} \Vert^2\Big)
\notag\\&\hskip8mm
+\frac{\rho_{i}}{2}\sum_{l}\Vert \pcG_{i}(X^{k+1,k+1}_{i,l})-\pcG_{i}(X^{k+1,k}_{i,l}) \Vert^2
+\big(\sigma^{k}_{i}+\rho_{i}\big)\Vert X^{k+1}_{i}-X^{k}_{i} \Vert^2
\Big].
\label{LPGSProximalXFOC}
\end{align}
Following procedures similar to the above, the updates
\eqref{LPGSProximalZl}, \eqref{LPGSProximalsvZpilA},
\eqref{LPGSProximalsvZnilA}, and 
\eqref{LPGSProximalYilA} yield
\begin{align}
-\sum_{i=1}^{N}\Big(
             \rho_{i}\big( 2 (Z^{k+1}_{l}-X^{k+1}_{i,l})
                          +\svZp^{k}_{i,l}-\svZn^{k}_{i,l}
                          \big)
            +\tau^{k}(Z^{k+1}_{l}-Z^{k}_{l})
            +\dvZp^{k}_{i,l}
            -\dvZn^{k}_{i,l}
        \Big)
\in
\partial_{Z_{l}} \Indicator_{[-u_{l},u_{l}]}(Z^{k+1}_{l}),
\label{LPGSProximalZkAFOCL2Norm}
\end{align}
\vspace{-5mm}
\begin{align} 
-(\svZpcoef^{k}_{i}+\rho_{i})\svZp^{k+1}_{i,l}
+\svZpcoef^{k}_{i}\svZp^{k}_{i,l}
-\rho_{i}(Z^{k+1}_{l}-X^{k+1}_{i,l})
-\dvZp^{k}_{i,l}
\in
 \partial_{\svZp_{i,l}}\Indicator_{[0, u_{\svZp_{i,l}}]}(\svZp^{k+1}_{i,l}),
\label{LPGSProximalsvZnilL2Norm}
\end{align}
\vspace{-5mm}
\begin{align} 
-(\svZncoef^{k}_{i}+\rho_{i})\svZn^{k+1}_{i,l}
+\svZncoef^{k}_{i}\svZn^{k}_{i,l}
-\rho_{i}(X^{k+1}_{i,l}-Z^{k+1}_{l})
-\dvZn^{k}_{i,l}
\in
 \partial_{\svZn_{i,l}}\Indicator_{[0, u_{\svZn_{i,l}}]}(\svZn^{k+1}_{i,l}),
\label{LPGSProximalsvZpilL2Norm}
\end{align}
\vspace{-5mm}
\begin{align} 
-(\gamma^{k}_{i}+\rho_{i}) Y^{k+1}_{i}
+\gamma^{k}_{i}Y^{k}_{i}
-\rho_{i} \pcG_{i}(X^{k+1}_{i})
-\dvG^{k}_{i}
\in
 \partial_{Y_{i}}\Indicator_{[0, u_{Y_{i}}]}(Y^{k+1}_{i}),   
\label{LPGSProximalYiAFOCL2Norm}
\end{align}
\vspace{-5mm}
\begin{align}
&
\sum_{i=1}^{N}\Big[
         \Indicator_{[-u_{},u_{}]}(Z^{k}_{})
        +\langle  \dvZp^{k}_{i}-\dvZn^{k}_{i}, Z^{k}_{}\rangle
        +\frac{\rho_{i}}{2}\Big(\Vert Z^{k}_{}-X^{k+1}_{i}+\svZp^{k}_{i} \Vert^2
                                +\Vert X^{k+1}_{i}-Z^{k}_{}+\svZn^{k}_{i} \Vert^2\Big)
                \Big]
\notag\\
\geq\,&
\sum_{i=1}^{N}\Big[
         \Indicator_{[-u_{},u_{}]}(Z^{k+1}_{})
        +\langle  \dvZp^{k}_{i}-\dvZn^{k}_{i}, Z^{k+1}_{}\rangle
        +\frac{\rho_{i}}{2}\Big(\Vert Z^{k+1}_{}-X^{k+1}_{i}+\svZp^{k}_{i} \Vert^2
                        +\Vert X^{k+1}_{i}-Z^{k+1}_{}+\svZn^{k}_{i} \Vert^2\Big)
                \Big]
\notag\\[-3pt]&
+\sum_{i=1}^{N}\big(\tau^{k}+\rho_{i}\big)\Vert Z^{k+1}_{}-Z^{k}_{} \Vert^2,
\label{LPGSProximalZFOC}
\end{align}
\vspace{-5mm}
\begin{align} 
&
 \sum_{i=1}^{N}\Big[
         \Indicator_{[0, u_{\svZp_{i}}]}(\svZp^{k}_{i}) 
        +\langle  \dvZp^{k}_{i}, \svZp^{k}_{i}\rangle
        +\frac{\rho_{i}}{2}\Vert Z^{k+1}_{}-X^{k+1}_{i}+\svZp^{k}_{i} \Vert^2
        \Big]
\notag\\[-3pt]
\geq\,&
 \sum_{i=1}^{N}\Big[
         \Indicator_{[0, u_{\svZp_{i}}]}(\svZp^{k+1}_{i}) 
        +\langle  \dvZp^{k}_{i}, \svZp^{k+1}_{i}\rangle
        +\frac{\rho_{i}}{2}\Vert Z^{k+1}_{}-X^{k+1}_{i}+\svZp^{k+1}_{i} \Vert^2
        \Big]
+ \sum_{i=1}^{N}\big(\svZpcoef^{k}_{i}+\frac{\rho_{i}}{2}\big) \Vert \svZp^{k+1}_{i}-\svZp^{k}_{i} \Vert^2,
\label{LPGSProximalsvZnkAFOCL2Norm}
\end{align}
\vspace{-3mm}
\begin{align} 
&
 \sum_{i=1}^{N}\Big[
         \Indicator_{[0, u_{\svZn_{i}}]}(\svZn^{k}_{i})
        +\langle  \dvZn^{k}_{i}, \svZn^{k}_{i}\rangle
        +\frac{\rho_{i}}{2}\Vert X^{k+1}_{i}-Z^{k+1}_{}+\svZn^{k}_{i} \Vert^2
\Big]
\notag\\[-3pt]
\geq\,&
 \sum_{i=1}^{N}\Big[
         \Indicator_{[0, u_{\svZn_{i}}]}(\svZn^{k+1}_{i}) 
        +\langle  \dvZn^{k}_{i}, \svZn^{k+1}_{i}\rangle
        +\frac{\rho_{i}}{2}\Vert X^{k+1}_{i}-Z^{k+1}_{}+\svZn^{k+1}_{i} \Vert^2
        \Big]
+\sum_{i=1}^{N}\big(\svZncoef^{k}_{i}+\frac{\rho_{i}}{2}\big)\Vert \svZn^{k+1}_{i}-\svZn^{k}_{i} \Vert^2,
\label{LPGSProximalsvZpkAFOCL2Norm}
\end{align}
and
\begin{align} 
&
 \sum_{i=1}^{N}\Big[ 
         \Indicator_{[0, u_{Y_{i}}]}(Y^{k}_{i}) 
        +\langle \dvG^{k}_{i}, Y^{k}_{i}\rangle
        +\frac{\rho_{i}}{2}\Vert \pcG_{i}(X^{k+1}_{i})+Y^{k}_{i} \Vert^2
        \Big]
\notag\\[-3pt]
\geq\,&
 \sum_{i=1}^{N}\Big[ 
         \Indicator_{[0, u_{Y_{i}}]}(Y^{k+1}_{i}) 
        +\langle \dvG^{k}_{i}, Y^{k+1}_{i}\rangle
        +\frac{\rho_{i}}{2}\Vert \pcG_{i}(X^{k+1}_{i})+Y^{k+1}_{i} \Vert^2
        \Big]
+\sum_{i=1}^{N}\big(\gamma^{k}_{i}+\frac{\rho_{i}}{2}\big) \Vert Y^{k+1}_{i}-Y^{k}_{i} \Vert^2.
\label{LPGSProximalZAFOCL2Norm}
\end{align}
Summation of \eqref{LPGSProximalXFOC}, \eqref{LPGSProximalZFOC}
through \eqref{LPGSProximalZAFOCL2Norm}
 and removal of the indicators satisfied by the primal updates give
\begin{lemma} 
The primal updates~\eqref{LPGSProximalXilA}, \eqref{LPGSProximalZl},  
\eqref{LPGSProximalsvZpilA}, \eqref{LPGSProximalsvZnilA}, 
and \eqref{LPGSProximalYilA} yield, for all $k$,
\begin{align*} 
&
 \sum_{i=1}^{N}\Big[
         f(X^{k}_{i})
        +\big\langle \dvZp^{k}_{i}, Z^{k}_{}-X^{k}_{i}+\svZp^{k}_{i} \big\rangle
        +\big\langle \dvZn^{k}_{i}, X^{k}_{i}-Z^{k}_{}+\svZn^{k}_{i} \big\rangle
        +\big\langle \dvG^{k}_{i}, \pcG_{i}(X^{k}_{i})+Y^{k}_{i} \big\rangle
\notag\\[-3pt]&\hskip8mm
        +\frac{\rho_{i}}{2}\Big(\Vert Z^{k}_{}-X^{k}_{i}+\svZp^{k}_{i} \Vert^2
                +\Vert X^{k}_{i}-Z^{k}_{}+\svZn^{k}_{i} \Vert^2
                +\Vert \pcG_{i}(X^{k}_{i})+Y^{k}_{i} \Vert^2\Big)
\Big]
\notag\\[-3pt]
\geq\,&
 \sum_{i=1}^{N}\Big[
        f(X^{k+1}_{i})
        +\big\langle  \dvZp^{k}_{i}, Z^{k+1}_{}-X^{k+1}_{i}+\svZp^{k+1}_{i}\big\rangle
        +\big\langle  \dvZn^{k}_{i}, X^{k+1}_{i}-Z^{k+1}_{}+\svZn^{k+1}_{i}\big\rangle
        +\big\langle \dvG^{k}_{i}, \pcG_{i}(X^{k+1}_{i})+Y^{k+1}_{i} \big\rangle
\notag\\[-3pt]&\hskip8mm
        +\frac{\rho_{i}}{2}\Big(\Vert Z^{k+1}_{}-X^{k+1}_{i}+\svZp^{k+1}_{i} \Vert^2
                        +\Vert X^{k+1}_{i}-Z^{k+1}_{}+\svZn^{k+1}_{i} \Vert^2
                        +\Vert \pcG_{i}(X^{k+1}_{i})+Y^{k+1}_{i} \Vert^2 \Big)
\Big]
\notag\\&
+\sum_{i=1}^{N}\Big[
                 \frac{\rho_{i}}{2}\sum_{l}\Vert \pcG_{i}(X^{k+1,k+1}_{i,l})-\pcG_{i}(X^{k+1,k}_{i,l}) \Vert^2
                +\big(\sigma^{k}_{i}+\rho_{i}\big)\Vert X^{k+1}_{i}-X^{k}_{i} \Vert^2
                +\big(\tau^{k}+\rho_{i}\big)\Vert Z^{k+1}_{}-Z^{k}_{} \Vert^2
                \notag\\&\hskip12mm
                +\big(\svZpcoef^{k}_{i}+\frac{\rho_{i}}{2}\big) \Vert \svZp^{k+1}_{i}-\svZp^{k}_{i} \Vert^2
                +\big(\svZncoef^{k}_{i}+\frac{\rho_{i}}{2}\big)\Vert \svZn^{k+1}_{i}-\svZn^{k}_{i} \Vert^2
                +\big(\gamma^{k}_{i}+\frac{\rho_{i}}{2}\big) \Vert Y^{k+1}_{i}-Y^{k}_{i} \Vert^2
\Big].
\end{align*}
\label{LPGSProximalFOCp}
\end{lemma}
This is the fundamental lemma for us to explore algorithmic modeling of the dual variable updates and 
analyze the initialization and convergence of the algorithm composed of \eqref{LPGSProximalXilA},
\eqref{LPGSProximalZl}, \eqref{LPGSProximalsvZpilA}, \eqref{LPGSProximalsvZnilA},
\eqref{LPGSProximalYilA}, and the dual updates to be introduced below.
Motivated by the mathematical structure of Lemma~\ref{LPGSProximalFOCp}
and satisfaction of the constraints in the consensus primal
problem~\eqref{LPGSPrimalProblemGlobalConsensusRing}
subject to the extended constraints of equality~\eqref{LPGSConstraintsithCBSlacked},
the following simple models are proposed to update the dual variables,
\begin{align}
&
\dvZp^{k+1}_{i}:=\dvZp^{k}_{i}-\alpha^{k+1}_{\dvZp_{i}}\big(Z^{k+1}_{}-X^{k+1}_{i}+\svZp^{k+1}_{i}\big),\ \
\dvZn^{k+1}_{i}:=\dvZn^{k}_{i}-\alpha^{k+1}_{\dvZn_{i}}\big(X^{k+1}_{i}-Z^{k+1}_{}+\svZn^{k+1}_{i}\big),
\notag\\&
\dvG^{k+1}_{i}:=\dvG^{k}_{i} - \alpha^{k+1}_{\dvG_{i}}\big(\pcG_{i}(X^{k+1}_{i})+Y^{k+1}_{i}\big),
\label{LPGSDualUpdateRules}
\end{align}
where the dual coefficients, 
$\alpha^{k+1}_{\dvZp_{i}}$, $\alpha^{k+1}_{\dvZn_{i}}$, and 
$\alpha^{k+1}_{\dvG_{i}}$
are positive scalar constants.
These models are equivalent to the proximal point algorithms of
\begin{align}
&
\dvZp^{k+1}_{i}
=\argmin_{\dvZp_{i}}\Big\{
  (2\alpha^{k+1}_{\dvZp_{i}})^{-1}\Vert \dvZp_{i}-\dvZp^{k}_{i}\Vert^2
 +\big\langle \dvZp_{i}, Z^{k+1}_{}-X^{k+1}_{i}+\svZp^{k+1}_{i}\big\rangle\Big\},
\notag\\[-0pt]&
\dvZn^{k+1}_{i}
=\argmin_{\dvZn_{i}}\Big\{
  (2\alpha^{k+1}_{\dvZn_{i}})^{-1}\Vert \dvZn_{i}-\dvZn^{k}_{i}\Vert^2
 +\big\langle \dvZn_{i}, X^{k+1}_{i}-Z^{k+1}_{}+\svZn^{k+1}_{i}\big\rangle\Big\},
\notag\\[-0pt]&
\dvG^{k+1}_{i}
=\argmin_{\dvG_{i}}\Big\{
  (2\alpha^{k+1}_{\dvG_{i}})^{-1}\Vert \dvG_{i}-\dvG^{k}_{i}\Vert^2
 +\big\langle \dvG_{i}, \pcG_{i}(X^{k+1}_{i})+Y^{k+1}_{i}\big\rangle\Big\}.
\label{LPSGProx4Duals}
\end{align}
These expressions indicate that the dual coefficients,
$\alpha^{k+1}_{\dvZp_{i}}$, $\alpha^{k+1}_{\dvZn_{i}}$, and
$\alpha^{k+1}_{\dvG_{i}}$
should be low 
and the constraints of $\pcG_{i}(Z)\leq 0$ should be preconditioned
via scaling-down if needed
such that the functions minimized are quite strongly convex,
in order to strengthen the boundedness of the dual sequences,
$\{\dvZp^{k}_{i}$, $\dvZn^{k}_{i}$, $\dvG^{k}_{i}\}_{k}$,
as to be analyzed in Sec.~\ref{sec:Initialization}.
To enhance possibly further the convexity,
the above algorithms may be extended by an introduction of direct
 couplings between $\dvZp_{i}-\dvZp^{k}_{i}$, $\dvZn_{i}-\dvZn^{k}_{i}$, 
and $\dvG_{i}-\dvG^{k}_{i}$
via a positively definite symmetric matrix of large eigenvalues; 
minimization of such a single quadratic convex function
is then employed to update the dual variables.
The above are dual descent rules for the augmented Lagrangian of linear programming,
in contrast to conventional dual ascent. 
Their necessity is to become clear by Lemma~\ref{LPGSPDp}
and its consequences to be inferred below. 

Substitution of \eqref{LPGSDualUpdateRules} into
Lemma~\ref{LPGSProximalFOCp} gives
\begin{lemma}
For all $k$,
\begin{align*} 
&
 \sum_{i=1}^{N}\Big[
         f(X^{k}_{i})
        +\big\langle \dvZp^{k}_{i}, Z^{k}_{}-X^{k}_{i}+\svZp^{k}_{i} \big\rangle
        +\big\langle \dvZn^{k}_{i}, X^{k}_{i}-Z^{k}_{}+\svZn^{k}_{i} \big\rangle
        +\big\langle \dvG^{k}_{i}, \pcG_{i}(X^{k}_{i})+Y^{k}_{i} \big\rangle
\notag\\[-3pt]&\hskip8mm
        +\frac{\rho_{i}}{2}\Big(\Vert Z^{k}_{}-X^{k}_{i}+\svZp^{k}_{i} \Vert^2
                +\Vert X^{k}_{i}-Z^{k}_{}+\svZn^{k}_{i} \Vert^2
                +\Vert \pcG_{i}(X^{k}_{i})+Y^{k}_{i} \Vert^2\Big)
\Big]
\notag\\[-3pt]
\geq\,&
 \sum_{i=1}^{N}\Big[
        f(X^{k+1}_{i})
        +\big\langle  \dvZp^{k+1}_{i}, Z^{k+1}_{}-X^{k+1}_{i}+\svZp^{k+1}_{i}\big\rangle
        +\big\langle  \dvZn^{k+1}_{i}, X^{k+1}_{i}-Z^{k+1}_{}+\svZn^{k+1}_{i}\big\rangle
        +\big\langle \dvG^{k+1}_{i}, \pcG_{i}(X^{k+1}_{i})+Y^{k+1}_{i} \big\rangle
\notag\\[-3pt]&\hskip8mm
        +\frac{\rho_{i}}{2}\Big(\Vert Z^{k+1}_{}-X^{k+1}_{i}+\svZp^{k+1}_{i} \Vert^2
                        +\Vert X^{k+1}_{i}-Z^{k+1}_{}+\svZn^{k+1}_{i} \Vert^2
                        +\Vert \pcG_{i}(X^{k+1}_{i})+Y^{k+1}_{i} \Vert^2 \Big)
\Big]
\notag\\[-3pt]&
+\sum_{i=1}^{N}\Big[
         \alpha^{k+1}_{\dvZp_{i}}\big\Vert Z^{k+1}_{}-X^{k+1}_{i}+\svZp^{k+1}_{i} \big\Vert^2
        +\alpha^{k+1}_{\dvZn_{i}}\big\Vert X^{k+1}_{i}-Z^{k+1}_{}+\svZn^{k+1}_{i} \big\Vert^2
        +\alpha^{k+1}_{\dvG_{i}}\big\Vert \pcG_{i}(X^{k+1}_{i})+Y^{k+1}_{i} \big\Vert^2
        \Big]
\notag\\[-3pt]&
+\sum_{i=1}^{N}\Big[
                 \frac{\rho_{i}}{2}\sum_{l}\Vert \pcG_{i}(X^{k+1,k+1}_{i,l})-\pcG_{i}(X^{k+1,k}_{i,l}) \Vert^2
                +\big(\sigma^{k}_{i}+\rho_{i}\big)\Vert X^{k+1}_{i}-X^{k}_{i} \Vert^2
                +\big(\tau^{k}+\rho_{i}\big)\Vert Z^{k+1}_{}-Z^{k}_{} \Vert^2
                \notag\\[-2pt]&\hskip12mm
                +\big(\svZpcoef^{k}_{i}+\frac{\rho_{i}}{2}\big) \Vert \svZp^{k+1}_{i}-\svZp^{k}_{i} \Vert^2
                +\big(\svZncoef^{k}_{i}+\frac{\rho_{i}}{2}\big)\Vert \svZn^{k+1}_{i}-\svZn^{k}_{i} \Vert^2
                +\big(\gamma^{k}_{i}+\frac{\rho_{i}}{2}\big) \Vert Y^{k+1}_{i}-Y^{k}_{i} \Vert^2
\Big].
\end{align*}
\label{LPGSPDp}
\end{lemma}

For fixed $K\in\mathbb{N}$,  $\sum_{k=0}^{K-1}$ is operated on
 Lemma~\ref{LPGSPDp}; the result is summarized in Lemma~\ref{LPGSPDSum2Kp}.
\begin{lemma} 
For all $K\in \mathbb{N}$, 
\begin{align*} 
&
 \sum_{i=1}^{N}
 \Big[f(X^{0}_{i})
        +\big\langle \dvZp^{0}_{i}, Z^{0}_{}-X^{0}_{i}+\svZp^{0}_{i} \big\rangle
        +\big\langle \dvZn^{0}_{i}, X^{0}_{i}-Z^{0}_{}+\svZn^{0}_{i} \big\rangle
        +\big\langle \dvG^{0}_{i}, \pcG_{i}(X^{0}_{i})+Y^{0}_{i} \big\rangle
        \notag\\[-3pt]&\hskip8mm
        +\frac{\rho_{i}}{2}\Big(\Vert Z^{0}_{}-X^{0}_{i}+\svZp^{0}_{i} \Vert^2
                +\Vert X^{0}_{i}-Z^{0}_{}+\svZn^{0}_{i} \Vert^2
                +\Vert \pcG_{i}(X^{0}_{i})+Y^{0}_{i} \Vert^2 \Big)
        \Big] 
\notag\\
\geq\,&
 \sum_{i=1}^{N}
 \Big[f(X^{K}_{i})
        +\big\langle \dvZp^{K}_{i}, Z^{K}_{}-X^{K}_{i}+\svZp^{K}_{i} \big\rangle
        +\big\langle \dvZn^{K}_{i}, X^{K}_{i}-Z^{K}_{}+\svZn^{K}_{i} \big\rangle
        +\big\langle \dvG^{K}_{i}, \pcG_{i}(X^{K}_{i})+Y^{K}_{i} \big\rangle
        \notag\\[-3pt]&\hskip8mm
        +\frac{\rho_{i}}{2}\Big(\Vert Z^{K}_{}-X^{K}_{i}+\svZp^{K}_{i} \Vert^2
                +\Vert X^{K}_{i}-Z^{K}_{}+\svZn^{K}_{i} \Vert^2
                +\Vert \pcG_{i}(X^{K}_{i})+Y^{K}_{i} \Vert^2 \Big)
        \Big]   
\notag\\[-0pt]&
+\sum_{i=1}^{N}\sum_{k=0}^{K-1}\Big[
         \alpha^{k+1}_{\dvZp_{i}}\big\Vert Z^{k+1}_{}-X^{k+1}_{i}+\svZp^{k+1}_{i} \big\Vert^2
        +\alpha^{k+1}_{\dvZn_{i}}\big\Vert X^{k+1}_{i}-Z^{k+1}_{}+\svZn^{k+1}_{i} \big\Vert^2
        +\alpha^{k+1}_{\dvG_{i}}\big\Vert \pcG_{i}(X^{k+1}_{i})+Y^{k+1}_{i} \big\Vert^2
        \Big]
\notag\\[-0pt]&
+\sum_{i=1}^{N}\sum_{k=0}^{K-1}\Big[
                 \frac{\rho_{i}}{2}\sum_{l}\Vert \pcG_{i}(X^{k+1,k+1}_{i,l})-\pcG_{i}(X^{k+1,k}_{i,l}) \Vert^2
                +\big(\sigma^{k}_{i}+\rho_{i}\big)\Vert X^{k+1}_{i}-X^{k}_{i} \Vert^2
                +\big(\tau^{k}+\rho_{i}\big)\Vert Z^{k+1}_{}-Z^{k}_{} \Vert^2
                \notag\\[-2pt]&\hskip19mm
                +\big(\svZpcoef^{k}_{i}+\frac{\rho_{i}}{2}\big) \Vert \svZp^{k+1}_{i}-\svZp^{k}_{i} \Vert^2
                +\big(\svZncoef^{k}_{i}+\frac{\rho_{i}}{2}\big)\Vert \svZn^{k+1}_{i}-\svZn^{k}_{i} \Vert^2
                +\big(\gamma^{k}_{i}+\frac{\rho_{i}}{2}\big) \Vert Y^{k+1}_{i}-Y^{k}_{i} \Vert^2
                \Big].
\end{align*}
\label{LPGSPDSum2Kp}
\end{lemma}

\begin{lemma} 
The sequences,
\begin{align*}
 \big\{ \,& X^{k}_{i}, Z^{k}, \svZp^{k}_{i}, \svZn^{k}_{i}, Y^{k}_{i},
        \dvZp^{k}_{i}, \dvZn^{k}_{i}, \dvG^{k}_{i},
          X^{k+1}_{i}, Z^{k+1}, \svZp^{k+1}_{i}, \svZn^{k+1}_{i}, Y^{k+1}_{i},
        \dvZp^{k+1}_{i}, \dvZn^{k+1}_{i}, \dvG^{k+1}_{i},\
    i=1,\ldots, N \big\}_{k}
\end{align*}
are supposedly bounded via the primal and dual update rules
and their initialization and parameter values,
there exist convergent subsequences, 
\begin{align*}
\big\{\,& X^{k_j}_{i}, Z^{k_j}, \svZp^{k_j}_{i}, \svZn^{k_j}_{i}, Y^{k_j}_{i},
        \dvZp^{k_j}_{i}, \dvZn^{k_j}_{i}, \dvG^{k_j}_{i},
        X^{k_j+1}_{i}, Z^{k_j+1}, \svZp^{k_j+1}_{i}, \svZn^{k_j+1}_{i}, Y^{k_j+1}_{i},
        \dvZp^{k_j+1}_{i}, \dvZn^{k_j+1}_{i}, \dvG^{k_j+1}_{i},
        \  i=1,\ldots, N \big\}_{j}
\end{align*}
 and accumulation points such that for all $i$,
\begin{align*}
&
\lim_{j\rightarrow\infty} X^{k_j}_{i}=X^{\infty}_{i},\
\lim_{j\rightarrow\infty} X^{k_j+1}_{i}=X^{\infty+}_{i},\
\lim_{j\rightarrow\infty} Z^{k_j}=Z^{\infty},\
\lim_{j\rightarrow\infty} Z^{k_j+1}=Z^{\infty+},\
\lim_{j\rightarrow\infty} \svZp^{k_j}_{i}=\svZp^{\infty}_{i},
\notag\\&
\lim_{j\rightarrow\infty} \svZp^{k_j+1}_{i}=\svZp^{\infty+}_{i},\
\lim_{j\rightarrow\infty} \svZn^{k_j}_{i}=\svZn^{\infty}_{i},\
\lim_{j\rightarrow\infty} \svZn^{k_j+1}_{i}=\svZn^{\infty+}_{i},\
\lim_{j\rightarrow\infty} Y^{k_j}_{i}=Y^{\infty}_{i},\
\lim_{j\rightarrow\infty} Y^{k_j+1}_{i}=Y^{\infty+}_{i},
\notag\\&
\lim_{j\rightarrow\infty} \dvZp^{k_j}_{i}=\dvZp^{\infty}_{i},\
\lim_{j\rightarrow\infty} \dvZp^{k_j+1}_{i}=\dvZp^{\infty+}_{i},\
\lim_{j\rightarrow\infty} \dvZn^{k_j}_{i}=\dvZn^{\infty}_{i},\
\lim_{j\rightarrow\infty} \dvZn^{k_j+1}_{i}=\dvZn^{\infty+}_{i},\
\lim_{j\rightarrow\infty} \dvG^{k_j}_{i}=\dvG^{\infty}_{i},\
\lim_{j\rightarrow\infty} \dvG^{k_j+1}_{i}=\dvG^{\infty+}_{i}.
\end{align*}
\label{LPGSPDAccumulationPointp}
\end{lemma}

\begin{lemma} 
Under $\rho_{i}> 0$, 
 $\{\sigma^{k}_{i}$, $\tau^{k}$, $\svZpcoef^{k}_{i}$, $\svZncoef^{k}_{i}$, $\gamma^{k}_{i}$, 
$\alpha^k_{\dvZp}$, $\alpha^{k}_{\dvZn_{i}}$, $\alpha^{k}_{\dvG_{i}}\}_{k}$
greater than some positive finite values, respectively,
 and the boundedness of the sequences supposed in Lemma~\ref{LPGSPDAccumulationPointp},
  we have
\begin{enumerate}
 \item[(i)] 
\vspace{-1mm}
$\forall i$,
\begin{align*}
\sum_{k=1}^{\infty}\big\Vert Z^{k}_{}-X^{k}_{i}+\svZp^{k}_{i} \big\Vert^2<\infty,\ \
\sum_{k=1}^{\infty}\big\Vert X^{k}_{i}-Z^{k}_{}+\svZn^{k}_{i} \big\Vert^2<\infty,\ \
\sum_{k=1}^{\infty}\big\Vert \pcG_{i}(X^{k}_{i})+Y^{k}_{i} \big\Vert^2<\infty;
\end{align*}
 \item[(ii)] 
\vspace{-2mm}
$\forall i$,
\begin{align*}
 &
 \sum_{k=0}^{\infty}\big\Vert X^{k+1}_{i}-X^{k}_{i} \big\Vert^2<\infty,\ \
 \sum_{k=0}^{\infty}\big\Vert Z^{k+1}_{}-Z^{k}_{} \big\Vert^2<\infty,\ \
 \sum_{k=0}^{\infty} \big\Vert \svZp^{k+1}_{i}-\svZp^{k}_{i} \big\Vert^2<\infty,\ \
 \sum_{k=0}^{\infty}\big\Vert \svZn^{k+1}_{i}-\svZn^{k}_{i} \big\Vert^2<\infty,
\notag\\[-1pt]&
 \sum_{k=0}^{\infty}\big\Vert Y^{k+1}_{i}-Y^{k}_{i} \big\Vert^2<\infty.
\end{align*}
\end{enumerate}
\label{LPGSPDSumBoundedp}
\end{lemma}
\begin{proof}
 Take $K=k_j$ in Lemma~\ref{LPGSPDSum2Kp} and
 apply Lemma~\ref{LPGSPDAccumulationPointp}.
\end{proof}

\begin{lemma} 
Under the conditions listed in Lemma~\ref{LPGSPDSumBoundedp}, we obtain the following:
\begin{enumerate}
 \item[(i)] 
\vspace{-2mm}
$\forall i$,
\begin{align*}
\lim_{k\rightarrow\infty} ( X^{k}_{i}-Z^{k})=0,\ \
\lim_{k\rightarrow\infty} \svZp^{k}_{i}=\lim_{k\rightarrow\infty} \svZn^{k}_{i}=0,\ \
\lim_{k\rightarrow\infty} \big( \pcG_{i}(X^{k}_{i})+Y^{k}_{i} \big)=0.
\end{align*}
\item[(ii)] 
\vspace{-2mm}
$\forall i$,
\begin{align*}
&
\lim_{k\rightarrow\infty} ( X^{k+1}_{i}-X^{k}_{i} )=0,\ \
\lim_{k\rightarrow\infty} ( Z^{k+1}-Z^{k} )=0,\ \
\lim_{k\rightarrow\infty} ( \svZp^{k+1}-\svZp^{k} )=0,\ \
\lim_{k\rightarrow\infty} ( \svZn^{k+1}-\svZn^{k} )=0,
\notag\\&
\lim_{k\rightarrow\infty} ( Y^{k+1}_{i}-Y^{k}_{i} )=0.
\end{align*}
\item[(iii)] 
\vspace{-2mm}
$\forall i$,
\begin{align*}
&
 X^{\infty+}_{i}=X^{\infty}_{i}=Z^{\infty}=Z^{\infty+},\ \
 \svZp^{\infty+}_{i}=\svZp^{\infty}_{i}
 =\svZn^{\infty+}_{i}=\svZn^{\infty}_{i}= 0,\ \
\pcG_{i}(X^{\infty}_{i})+Y^{\infty}_{i}=0,\ \
 Y^{\infty+}_{i}=Y^{\infty}_{i}\geq 0.
\end{align*}
The accumulation point $Z^{\infty}$ is a feasible solution of 
the primal problem~\eqref{LPGSPrimalProblemOriginal}.
\item[(iv)] 
The stationarity conditions are
\vspace{-1mm}
\begin{align*} 
&
-(f_{l})^T
+\dvZp^{\infty}_{i,l}
-\dvZn^{\infty}_{i,l}
-(\pcG_{i,l})^T\dvG^{\infty}_{i}
\in \partial_{X_{i,l}}\Indicator_{[-u_{l}, u_{l}]}(X^{\infty}_{i,l}),\ \
-\sum_{i=1}^{N}\big(\dvZp^{\infty}_{i}-\dvZn^{\infty}_{i}\big)
\in \partial_{Z_{}} \Indicator_{[-u_{},u_{}]}(Z^{\infty}_{}),
\notag\\&
-\dvZp^{\infty}_{i}
\in \partial_{\svZp_{i}}\Indicator_{[0, u_{\svZp_{i}}]}(0),\ \
-\dvZn^{\infty}_{i}
\in \partial_{\svZn_{i}}\Indicator_{[0, u_{\svZn_{i}}]}(0),\ \
-\dvG^{\infty}_{i}
\in \partial_{Y_{i}}\Indicator_{[0, u_{Y_{i}}]}(Y^{\infty}_{i}),
\end{align*}
where $f(X_{i})=\sum_{l=1}^{M}f_{l}X_{i,l}$ is used.
Next, 
\begin{align*}
\dvZp^{\infty}_{i}\geq 0,\ \ 
\dvZn^{\infty}_{i}\geq 0,\ \ 
\dvG^{\infty}_{i}\geq 0,\ \ \big\langle \dvG^{\infty}_{i}, \pcG_{i}(X^{\infty}_{i})\big\rangle=0.
\end{align*} 
\end{enumerate}
\label{LPGSPDDeltaLimitsp}
\end{lemma}
\begin{proof}
 (i), (ii), and (iii) can be verified directly;
 $Y^{\infty}_{i}\geq 0$ is guaranteed by the update rule~\eqref{LPGSProximalYilA}.
 The first part of (iv) comes from \eqref{LPGSProximalXilAFOC},
 \eqref{LPGSProximalZkAFOCL2Norm} through \eqref{LPGSProximalYiAFOCL2Norm},
  Lemma~\ref{LPGSPDAccumulationPointp}, (i), (ii), and (iii).
 The second part of (iv) follows from the first part of (iv)
 and how the upper bounds for the slack variables are set.
  \end{proof}
 
Motivated by Lemma~\ref{LPGSPDp}, the following Lagrangian function sequence is introduced.
\begin{definition}
For all $k$,
\begin{align*}
L^{k}:=\,&
\sum_{i=1}^{N}
 \Big[f(X^{k}_{i})
        +\big\langle \dvZp^{k}_{i}, Z^{k}_{}-X^{k}_{i}+\svZp^{k}_{i} \big\rangle
        +\big\langle \dvZn^{k}_{i}, X^{k}_{i}-Z^{k}_{}+\svZn^{k}_{i} \big\rangle
        +\big\langle \dvG^{k}_{i}, \pcG_{i}(X^{k}_{i})+Y^{k}_{i} \big\rangle
        \notag\\[-3pt]&\hskip8mm
        +\frac{\rho_{i}}{2}\Big(\Vert Z^{k}_{}-X^{k}_{i}+\svZp^{k}_{i} \Vert^2
                +\Vert X^{k}_{i}-Z^{k}_{}+\svZn^{k}_{i} \Vert^2
                +\Vert \pcG_{i}(X^{k}_{i})+Y^{k}_{i} \Vert^2\Big) \Big].
\end{align*}
\label{LPGSNonIncreasingFunctionDefn}
\end{definition}
$L^{k}$ has the same formal structure as that of the augmented Lagrangian $L$ 
defined through \eqref{LPGSLagranginDualConsensus} 
and \eqref{LPGSLagranginFunctionithBC}.
Lemma~\ref{LPGSNonIncreasingFunctionp} below follows from 
Lemmas~\ref{LPGSPDp}, \ref{LPGSPDAccumulationPointp},
 and \ref{LPGSPDDeltaLimitsp}.
\begin{lemma}
  $\forall k$, $L^{k}\geq L^{k+1}$.
 $\lim_{k\rightarrow\infty}L^{k}= N f(Z^{\infty})$,
$\lim_{k\rightarrow\infty} f(Z^{k})=f(Z^{\infty})$,
 $\lim_{k\rightarrow\infty} f(X^{k}_{i})=f(X^{\infty}_{i})=f(Z^{\infty})$.
\label{LPGSNonIncreasingFunctionp}
\end{lemma}
Lemma~\ref{LPGSNonIncreasingFunctionp} indicates that
the primal and dual sequences need to be initialized adequately in order to make the
algorithm feasible.

\begin{lemma}
For all $k$,
\begin{align*} 
 \sum_{i=1}^{N}f(X^{\ast}_{i})
\geq\,&
 \sum_{i=1}^{N}\s\Big[
         f(X^{k+1}_{i})
        +\big\langle \dvZp^{k}_{i}, Z^{k+1}_{}-X^{k+1}_{i}+\svZp^{k+1}_{i}\big\rangle
        +\big\langle \dvZn^{k}_{i}, X^{k+1}_{i}-Z^{k+1}_{}+\svZn^{k+1}_{i}\big\rangle
        +\big\langle \dvG^{k}_{i},  \pcG_{i}(X^{k+1}_{i})+Y^{k+1}_{i}\big\rangle
        \notag\\[-5pt]&\hskip7.4mm
        +\rho_{i}\Big( 
                 \big\langle Z^{k+1}_{}-X^{k+1}_{i}+\svZp^{k+1}_{i}, 
                                \svZp^{k+1}_{i} \big\rangle
                +\big\langle X^{k+1}_{i}-Z^{k+1}_{}+\svZn^{k+1}_{i}, 
                                \svZn^{k+1}_{i} \big\rangle
                \notag\\&\hskip17.4mm
                +\big\langle Z^{k+1}_{}-X^{k+1}_{i}+\svZp^{k}_{i}, 
                                        Z^{k+1}_{}-Z^{\ast}_{} \big\rangle
                -\big\langle Z^{k}_{}-X^{k+1}_{i}+\svZp^{k}_{i}, 
                                        X^{k+1}_{i}-X^{\ast}_{i} \big\rangle
                \notag\\&\hskip17.4mm
                +\big\langle X^{k+1}_{i}-Z^{k}_{}+\svZn^{k}_{i}, 
                                        X^{k+1}_{i}-X^{\ast}_{i} \big\rangle
                -\big\langle X^{k+1}_{i}-Z^{k+1}_{}+\svZn^{k}_{i}, 
                                        Z^{k+1}_{}-Z^{\ast}_{} \big\rangle
                \notag\\&\hskip17.4mm
                +\big\langle \pcG_{i}(X^{k+1}_{i})+Y^{k+1}_{i}, Y^{k+1}_{i}-Y^{\ast}_{i} \big\rangle
                +\sum_{l}\s\big\langle \pcG_{i}(X^{k+1,k+1}_{i,l})+Y^{k}_{i}, 
                                        \pcG_{i}(X^{k+1,k+1}_{i,l})-\pcG_{i}(X^{k+1,\ast}_{i,l}) \big\rangle
                                        \s\Big)
        \notag\\[-3pt]&\hskip7.4mm
        +\sigma^{k}_{i}\big\langle X^{k+1}_{i}-X^{k}_{i}, X^{k+1}_{i}-X^{\ast}_{i} \big\rangle
        +\tau^{k}\big\langle Z^{k+1}_{}-Z^{k}_{}, 
                                Z^{k+1}_{}-Z^{\ast}_{} \big\rangle
        +\svZpcoef^{k}_{i}\big\langle \svZp^{k+1}_{i}-\svZp^{k}_{i}, 
                        \svZp^{k+1}_{i} \big\rangle
        \notag\\&\hskip7.4mm
        +\svZncoef^{k}_{i}\big\langle \svZn^{k+1}_{i}-\svZn^{k}_{i}, 
                        \svZn^{k+1}_{i} \big\rangle 
        +\gamma^{k}_{i}\big\langle Y^{k+1}_{i}-Y^{k}_{i}, Y^{k+1}_{i}-Y^{\ast}_{i} \big\rangle
                        \Big].
\end{align*}
\label{LPGSPDastkp} 
\end{lemma}
\begin{proof}
 Find the first-order characteristics of 
\eqref{LPGSProximalXilA},
\eqref{LPGSProximalZl}, \eqref{LPGSProximalsvZpilA}, \eqref{LPGSProximalsvZnilA}, and
\eqref{LPGSProximalYilA}
under $\{X_{i,l}=X^{\ast}_{i,l}$, $Z=Z^{\ast}$, 
 $\svZp_{i}=\svZp^{\ast}_{i}=-(Z^{\ast}-X^{\ast}_{i})$, 
 $\svZn_{i}=\svZn^{\ast}_{i}=-(X^{\ast}_{i}-Z^{\ast})$,
$Y_{i}=Y^{\ast}_{i}=-\pcG_{i}(X^{\ast}_{i})\}$,
 sum them, and apply \eqref{LPGSKKTConditions}.
\end{proof}
\begin{lemma}
 $f(Z^{\infty})=f(Z^{\ast})$.
 $Z^{\infty}$ is an optimal solution of 
the primal problem~\eqref{LPGSPrimalProblemOriginal}.
\label{LPGSastinftyp}
\end{lemma}
\begin{proof}
Apply Lemmas~\ref{LPGSPDDeltaLimitsp}, \ref{LPGSNonIncreasingFunctionp},
and \ref{LPGSPDastkp}. 
\end{proof}

Summarizing the above lemmas, we have
\begin{theorem}
Suppose that
 $\rho_{i}> 0$,
  $\{\sigma^{k}_{i}$, $\tau^{k}$, $\svZpcoef^{k}_{i}$, $\svZncoef^{k}_{i}$, $\gamma^{k}_{i}$, 
$\alpha^k_{\dvZp}$, $\alpha^{k}_{\dvZn_{i}}$, $\alpha^{k}_{\dvG_{i}}\}_{k}$
are greater than some positive finite values, respectively,
the algorithm is feasible,
 and the primal and dual sequences are bounded.
 Then,
\begin{enumerate}
\item[(i)] 
$\big\{X^{k}_{i}$, $Z^{k}$, $\svZp^{k}$, $\svZn^{k}$, $Y^{k}_{i}\big\}_{k}$ 
converge collectively
to $\big\{X^{\infty}_{i}-Z^{\infty}=0$,
$\svZp^{\infty}_{i}=\svZn^{\infty}_{i}=0$,
$\pcG_{i}(X^{\infty}_{i})+Y^{\infty}_{i}=0$, 
$Y^{\infty}_{i}\geq 0$,
$i=1,\ldots, N\big\}$;
in particular, $\pcH_{i}(X^{\infty}_{i})=0$;
\item[(ii)] 
$\dvZp^{\infty}_{i}\geq 0$, $\dvZn^{\infty}_{i}\geq 0$;
$\dvG^{\infty}_{i}\geq 0$, 
$\big\langle \dvG^{\infty}_{i}, \pcG_{i}(X^{\infty}_{i})\big\rangle=0$,
$-(f_{l})^T
+\dvZp^{\infty}_{i,l}
-\dvZn^{\infty}_{i,l}
-(\pcG_{i,l})^T\dvG^{\infty}_{i}
\in \partial_{X_{i,l}}\Indicator_{[-u_{l}, u_{l}]}(X^{\infty}_{i,l})$,
$-\sum_{i=1}^{N}\big(\dvZp^{\infty}_{i}-\dvZn^{\infty}_{i}\big)
\in \partial_{Z_{}} \Indicator_{[-u_{},u_{}]}(Z^{\infty}_{})$;
\item[(iii)] 
$\lim_{k\rightarrow\infty} f(X^{k}_{i})=f(X^{\infty}_{i})$, 
$\lim_{k\rightarrow\infty} f(Z^{k})=f(Z^{\infty})$, and 
 $Z^{\infty}$ is an optimal solution of 
the primal problem~\eqref{LPGSPrimalProblemOriginal};
\item[(iv)] 
$\big\{L^{k}\big\}_{k}$ is decreasing and
$\frac{1}{N}\lim_{k\rightarrow\infty}  L^{k}=f(Z^{\infty})$.
\end{enumerate}
\label{LPGSConvergenceT}
\end{theorem}
Here, we use the phrase ``converge collectively to'',
because the results are proved on the basis of the convergent subsequences
of Lemma~\ref{LPGSPDAccumulationPointp} and
the existence of the
separate limits of $\big\{X^{k}_{i}$, $Z^{k}$, $\svZp^{k}$, $\svZn^{k}$, $Y^{k}_{i}\big\}_{k}$,
such as $\lim_{k\rightarrow\infty}X^{k}_{i}=X^{\infty}_{i}$ 
and $\lim_{k\rightarrow\infty}Z^{k}=Z^{\infty}$, 
is not yet established.
Considering that the update rules,
\eqref{LPGSProximalXilA}, \eqref{LPGSProximalZl},  
\eqref{LPGSProximalsvZpilA}, \eqref{LPGSProximalsvZnilA},
\eqref{LPGSProximalYilA},
and \eqref{LPGSDualUpdateRules}
are well-defined, provided that the initialization and parameters are fixed adequately,
the separate limits of 
$\big\{X^{k}_{i}$, $Z^{k}$, $\svZp^{k}$, $\svZn^{k}$, $Y^{k}_{i}$,
$\dvZp^{k}_{i}$, $\dvZn^{k}_{i}$, $\dvG^{k}_{i}\big\}_{k}$
are expected to exist.

It follows from Theorem~\ref{LPGSConvergenceT}
that $\{L^{k}\}_{k}$ is decreasing and converges to $N f(Z^{\infty})$. 
It is thus natural to identify the rate of  convergence of $\{L^{k}\}_{k}$ as the rate of convergence 
of the algorithm composed of 
\eqref{LPGSProximalXilA}, \eqref{LPGSProximalZl},  
\eqref{LPGSProximalsvZpilA}, \eqref{LPGSProximalsvZnilA},
\eqref{LPGSProximalYilA},
and \eqref{LPGSDualUpdateRules}.
However, the presence of the penalty terms poses a challenge to the estimate of the rate.
 Alternatively, it is observed that $X^{\infty}_{i}-Z^{\infty}=0$ is a direct consequence
of 
$\sum_{k=1}^{\infty}\Vert Z^{k}_{}-X^{k}_{i}+\svZp^{k}_{i} \Vert^2<\infty$ and
$\sum_{k=1}^{\infty}\Vert X^{k}_{i}-Z^{k}_{}+\svZn^{k}_{i} \Vert^2<\infty$
from Lemma~\ref{LPGSPDSumBoundedp}(i):
$\lim_{k\rightarrow\infty}(Z^{k}_{}-X^{k}_{i}+\svZp^{k}_{i})=Z^{\infty}-X^{\infty}_{i}$,
$\lim_{k\rightarrow\infty}(X^{k}_{i}-Z^{k}_{}+\svZn^{k}_{i})=X^{\infty}_{i}-Z^{\infty}$,
and thus, these bounded sums may be used to estimate the rate of convergence of the algorithm.
To this end, recall the well-known result of
$\sum_{k=1}^{\infty}k^{-q}<\infty$ if and only if $q>1$,
we may then infer that there exists $\Lambda>0$ such that
 $\Vert Z^{k}_{}-X^{k}_{i}+\svZp^{k}_{i} \Vert^2<\Lambda^2/k$
and $\Vert X^{k}_{i}-Z^{k}_{}+\svZn^{k}_{i} \Vert^2<\Lambda^2/k$ for almost all large $k\in\mathbb{N}$.
Consequently, $O(1/k^{1/2})$ is an estimate of the rate of convergence of the algorithm,
stated in Theorem~\ref{LPGSConvgRate}.
\begin{theorem}
(\textbf{$O(1/k^{1/2})$ rate of convergence}).
Suppose that
 $\rho_{i}>0$, 
 $\{\sigma^{k}_{i}$, $\tau^{k}$, $\svZpcoef^{k}_{i}$, $\svZncoef^{k}_{i}$, $\gamma^{k}_{i}$,
$\alpha^k_{\dvZp}$, $\alpha^{k}_{\dvZn_{i}}$, $\alpha^{k}_{\dvG_{i}}\}_{k}$
are greater than some positive finite values, respectively,
the algorithm is feasible,
 and the primal and dual sequences are bounded.
Then, $O(1/k^{1/2})$ is an estimate of the rate of convergence of the algorithm composed of 
\eqref{LPGSProximalXilA}, \eqref{LPGSProximalZl},  
\eqref{LPGSProximalsvZpilA}, \eqref{LPGSProximalsvZnilA},
\eqref{LPGSProximalYilA},
and \eqref{LPGSDualUpdateRules}.
\label{LPGSConvgRate} 
\end{theorem}
This estimate is rather crude;
a faster rate of convergence is expected and
a more rigorous argument is to be explored.
 
\section{Initialization, Parameter Estimates, Feasibility, and Boundedness}
\label{sec:Initialization}

We have the algorithm composed of the update rules~\eqref{LPGSProximalXilA},
\eqref{LPGSProximalZl},  
\eqref{LPGSProximalsvZpilA}, \eqref{LPGSProximalsvZnilA},
\eqref{LPGSProximalYilA},
and \eqref{LPGSDualUpdateRules}
whose convergence is argued under
its feasibility and the boundedness of the primal and dual sequences supposed.
The boundedness of the primal sequences,
$\{X^{k}_{i}$, $Z^{k}$, $\svZp^{k}_{i}$, $\svZn^{k}_{i}$,
$Y^{k}_{i}\}_{k}$ is guaranteed explicitly through their update rules.
Concerning the boundedness of the dual sequences, 
$\{\dvZp^{k}_{i}$, $\dvZn^{k}_{i}$, $\dvG^{k}_{i}\}_{k}$,
we intend to achieve it through adequate initialization 
and control parameter values, a strategy of descent.

To analyze the issue of feasibility via initialization and parameter values, 
we obtain from 
Lemma~\ref{LPGSPDSum2Kp} and Theorem~\ref{LPGSConvergenceT},
\begin{align} 
 L^{0}
\geq\,&
 Nf(Z^{\infty})
+\sum_{i=1}^{N}\sum_{k=1}^{\infty}\Big[
         \alpha^{k}_{\dvZp_{i}}\big\Vert Z^{k}_{}-X^{k}_{i}+\svZp^{k}_{i} \big\Vert^2
        +\alpha^{k}_{\dvZn_{i}}\big\Vert X^{k}_{i}-Z^{k}_{}+\svZn^{k}_{i} \big\Vert^2
        +\alpha^{k}_{\dvG_{i}}\big\Vert \pcG_{i}(X^{k}_{i})+Y^{k}_{i} \big\Vert^2
        \Big]
\notag\\[-3pt]&
+\sum_{i=1}^{N}\sum_{k=0}^{\infty}\Big[
                 \frac{\rho_{i}}{2}\sum_{l}\Vert \pcG_{i}(X^{k+1,k+1}_{i,l})-\pcG_{i}(X^{k+1,k}_{i,l}) \Vert^2
                +\big(\sigma^{k}_{i}+\rho_{i}\big)\Vert X^{k+1}_{i}-X^{k}_{i} \Vert^2
                +\big(\tau^{k}+\rho_{i}\big)\Vert Z^{k+1}_{}-Z^{k}_{} \Vert^2
                \notag\\[-3pt]&\hskip18mm
                +\big(\svZpcoef^{k}_{i}+\frac{\rho_{i}}{2}\big) \Vert \svZp^{k+1}_{i}-\svZp^{k}_{i} \Vert^2
                +\big(\svZncoef^{k}_{i}+\frac{\rho_{i}}{2}\big)\Vert \svZn^{k+1}_{i}-\svZn^{k}_{i} \Vert^2
                +\big(\gamma^{k}_{i}+\frac{\rho_{i}}{2}\big) \Vert Y^{k+1}_{i}-Y^{k}_{i} \Vert^2
                \Big],
\label{LPGSProximalXZYAFOCSumInfty}
\end{align}
where
\begin{align} 
L^{0}:=
 \sum_{i=1}^{N}
 \Big[\,&f(X^{0}_{i})
        +\big\langle \dvZp^{0}_{i}, Z^{0}_{}-X^{0}_{i}+\svZp^{0}_{i} \big\rangle
        +\big\langle \dvZn^{0}_{i}, X^{0}_{i}-Z^{0}_{}+\svZn^{0}_{i} \big\rangle
        +\big\langle \dvG^{0}_{i}, \pcG_{i}(X^{0}_{i})+Y^{0}_{i} \big\rangle
        \notag\\[-3pt]&
        +\frac{\rho_{i}}{2}\Big(\Vert Z^{0}_{}-X^{0}_{i}+\svZp^{0}_{i} \Vert^2
                +\Vert X^{0}_{i}-Z^{0}_{}+\svZn^{0}_{i} \Vert^2
                +\Vert \pcG_{i}(X^{0}_{i})+Y^{0}_{i} \Vert^2 \Big)
        \Big]. 
\label{LPGSProximalXZYAFOCSumL0}
\end{align}
It follows that the primal and dual sequences need to be initialized such that $L^{0}$ is 
great to make \eqref{LPGSProximalXZYAFOCSumInfty} feasible.
To this end, the following assignments are taken:
(a) 
$Z^{0}_{l,j}=\lambda_{Z}\, \text{sign}(f_{l,j})\, u_{l,j}$
where $\text{sign}(f_{l,j})=1$ when $f_{l,j}= 0$ and some components of $f_{l}$ are positive,
 $\lambda_{Z}=0.8$, say, to make $Z^{0}$ well in the interior of $[-u,u]$.
Then, 
$f(Z^{0})=\lambda_{Z}\sum_{l,j}\vert f_{l,j}\vert u_{l,j}$;
 $X^{0}_{i}=Z^{0}$ $\forall i$.
(b)
$\svZn^{0}_{i}=\svZp^{0}_{i}=u_{\dvZp_{i}}=2u+\epsilon_{Z}$ for certain $\epsilon_{Z}> 0$.
(c)
$Y^{0}_{i}=u_{Y_{i}}=\sum_{l=1}^{M}\sum_{j}\vert \pcG_{i,l,j}\vert\, u_{l,j}+\pcG_{i,0}
+\epsilon_{\pcG_{i}}$ for certain $\epsilon_{\pcG_{i}}> 0$. 
(d) Based on the structure of $L^{0}$,
the initial values of the dual sequences are taken as
\begin{align}
\dvZp^{0}_{i}=\lambda_{\dvZp_{i}}(Z^{0}-X^{0}_{i}+\svZp^{0}_{i}),\ \
\dvZn^{0}_{i}=\lambda_{\dvZn_{i}}(X^{0}_{i}-Z^{0}+\svZn^{0}_{i}),\ \
\dvG^{0}_{i}=\lambda_{\dvG_{i}}(\pcG_{i}(X^{0}_{i})+Y^{0}_{i});\ \ 
\{\lambda_{\dvZp_{i}},\lambda_{\dvZn_{i}},\lambda_{\dvG_{i}}\}\propto \rho_{i}.
\label{LPGSInitializationDual}
\end{align}
Here, the proportional coefficients are set 
proportional to $\rho_{i}$,
 motivated by the compatibility numerically among
 the quantities involved in the primal update rules~\eqref{LPGSProximalXilA},
\eqref{LPGSProximalZl},  \eqref{LPGSProximalsvZpilA}, \eqref{LPGSProximalsvZnilA},
and \eqref{LPGSProximalYilA}, and they should be large to produce great $L^{0}$.
With these specific choices, $L^{0}$ reduces to
\begin{align}
 L^{0}=\sum_{i=1}^{N}\Big[
                 f(Z^{0}_{})
                +\big(\lambda_{\dvZp_{i}}+\frac{\rho_{i}}{2}\big)\Vert \svZp^{0}_{i} \Vert^2
                +\big(\lambda_{\dvZn_{i}}+\frac{\rho_{i}}{2}\big)\Vert \svZn^{0}_{i} \Vert^2
                +\big(\lambda_{\dvG_{i}} +\frac{\rho_{i}}{2}\big)\Vert \pcG_{i}(Z^{0}_{})+Y^{0}_{i} \Vert^2
\Big].
\label{LPGSL0}
\end{align}
This expression illustrates one important role played by the slack variables,
especially
$\{\svZp_{i}$,\,$\svZn_{i}\}$
introduced via the extended constraints of equality in that 
the choices of $\{\svZp^{0}_{i}$,\,$\svZn^{0}_{i}, Y^{0}_{i}\}$ can make $L^{0}$ great
under $X^{0}_{i}=Z^{0}$.

Next,
both \eqref{LPGSProximalXZYAFOCSumInfty} and \eqref{LPGSInitializationDual}
indicate that
\begin{align}
\lambda_{\dvZp_{i}}\gg \alpha^{k}_{\dvZp_{i}},\ \
\lambda_{\dvZn_{i}}\gg \alpha^{k}_{\dvZn_{i}},\ \ 
\lambda_{\dvG_{i}} \gg \alpha^{k}_{\dvG_{i}}.
\label{LPGSDualControlParametersEstimated}
\end{align}

The numerical values of $\pcG_{i,l}$
and the proximal parameters
 $\{\sigma^{k}_{i}$,\,$\tau^{k}$,\,$\svZpcoef^{k}_{i}$,  
 $\svZncoef^{k}_{i}$,\,$\gamma^{k}_{i}\}$ control directly
the magnitude of
\vspace{-2mm}
\begin{align}
&
\sum_{i=1}^{N}\sum_{k=0}^{\infty}\Big[
                 \frac{\rho_{i}}{2}\sum_{l}\Vert \pcG_{i}(X^{k+1,k+1}_{i,l})-\pcG_{i}(X^{k+1,k}_{i,l}) \Vert^2
                +\big(\sigma^{k}_{i}+\rho_{i}\big)\Vert X^{k+1}_{i}-X^{k}_{i} \Vert^2
                +\big(\tau^{k}+\rho_{i}\big)\Vert Z^{k+1}_{}-Z^{k}_{} \Vert^2
                \notag\\[-3pt]&\hskip14mm
                +\big(\svZpcoef^{k}_{i}+\frac{\rho_{i}}{2}\big) \Vert \svZp^{k+1}_{i}-\svZp^{k}_{i} \Vert^2
                +\big(\svZncoef^{k}_{i}+\frac{\rho_{i}}{2}\big)\Vert \svZn^{k+1}_{i}-\svZn^{k}_{i} \Vert^2
                +\big(\gamma^{k}_{i}+\frac{\rho_{i}}{2}\big) \Vert Y^{k+1}_{i}-Y^{k}_{i} \Vert^2
                \Big]
\label{LPGSProximalXZYAFOCSum2nd}
\end{align}
in \eqref{LPGSProximalXZYAFOCSumInfty}.
It is clear that the constraints $\{\pcG_{i}(Z)\leq 0\}$
should be preconditioned via scaling-down, if required, to control the magnitude 
of the first term in \eqref{LPGSProximalXZYAFOCSum2nd} involving the sums of
$\Vert \pcG_{i,l}(X^{k+1}_{i,l}-X^{k}_{i,l})\Vert^2$.
To estimate the effect of the proximal parameter values and infer the consequences, 
we resort to the special solutions of   
\eqref{LPGSProximalXilA},
\eqref{LPGSProximalZl},  \eqref{LPGSProximalsvZpilA}, \eqref{LPGSProximalsvZnilA},
and \eqref{LPGSProximalYilA},
\vspace{-2mm}
\begin{align} 
(\sigma^{k}_{i}+2\rho_{i})(X^{k+1}_{i,l}-X^{k}_{i,l})
=\,&
-\Big[ {1}_{m_{l}\times m_{l}}
       -\frac{\rho_{i}(\pcG_{i,l})^T\pcG_{i,l}}{\sigma^{k}_{i}+2\rho_{i}}
        \Big]
   \Big[
         2\rho_{i}(X^{k}_{i,l}-Z^{k}_{l})
        +\rho_{i}(\svZn^{k}_{i,l}-\svZp^{k}_{i,l})
        \notag\\[-0pt]&\hskip25mm
        +\dvZn^{k}_{i,l}-\dvZp^{k}_{i,l}
        +(\pcG_{i,l})^T\dvG^{k}_{i}
        +(f_{l})^T
        +\rho_{i}(\pcG_{i,l})^T(\pcG_{i}(X^{k+1,k}_{i,l})+Y^{k}_{i})
        \Big],
\label{LPGSdeltaXkiSoluIterated}
\end{align}
which is a second-order approximate solution adequate under 
${1}_{m_{l}\times m_{l}} \gg \rho_{i}/(\sigma^{k}_{i}+2\rho_{i})(\pcG_{i,l})^T\pcG_{i,l}$,
\vspace{-3mm}
\begin{align}
\Big(\sum_{i=1}^{N}(\tau^{k}+2\rho_{i})\Big)(Z^{k+1}_{l}-Z^{k}_{l})
=\,&
 \sum_{i=1}^{N}\frac{1}{\sigma^{k}_{i}+2\rho_{i}}
    \bigg\{
     \sigma^{k}_{i}
              \Big[      2\rho_{i}\big(X^{k}_{i,l}-Z^{k}_{l}\big)
                        +\rho_{i}\big(\svZn^{k}_{i,l}-\svZp^{k}_{i,l}\big)
                        +\dvZn^{k}_{i,l}
                        -\dvZp^{k}_{i,l}\Big]
\notag\\[-3pt]&\hskip23mm
    +2\rho_{i}
        \Big[
                -(f_{l})^T
                -(\pcG_{i,l})^T\dvG^{k}_{i}
                -\rho_{i}(\pcG_{i,l})^T\big(\pcG_{i}(X^{k+1,k}_{i,l})+Y^{k}_{i}\big)
                \Big]
        \bigg\},
\label{LPGSdeltaZkSoluIterated}
\end{align}
\vspace{-5mm}
\begin{align} 
\svZp^{k+1}_{i}-\svZp^{k}_{i}
=
-\frac{\rho_{i}\svZp^{k}_{i}+\rho_{i}(Z^{k+1}_{}-X^{k+1}_{i})+\dvZp^{k}_{i}}
      {\svZpcoef^{k}_{i}+\rho_{i}},
\label{LPGSdeltasvZnkSoluIterated}
\end{align}
\vspace{-5mm}
\begin{align} 
 \svZn^{k+1}_{i}-\svZn^{k}_{i}
=
-\frac{\rho_{i}\svZn^{k}_{i}+\rho_{i}(X^{k+1}_{i}-Z^{k+1}_{})+\dvZn^{k}_{i}}
      {\svZncoef^{k}_{i}+\rho_{i}},
\label{LPGSdeltasvZpkSoluIterated}
\end{align}
and
\vspace{-4mm}
\begin{align} 
Y^{k+1}_{i}-Y^{k}_{i}
=
-\frac{\rho_{i} Y^{k}_{i}+\rho_{i} \pcG_{i}(X^{k+1}_{i})+\dvG^{k}_{i}}
      {\gamma^{k}_{i}+\rho_{i}}.
\label{LPGSdeltasvGkSoluIterated}
\end{align}

Let us consider the consequence of \eqref{LPGSdeltaXkiSoluIterated},
\begin{align}
&
\sum_{i=1}^{N}\sum_{k=0}^{\infty}
        \big(\sigma^{k}_{i}+\rho_{i}\big)\Vert X^{k+1}_{i}-X^{k}_{i} \Vert^2
\notag\\
\leq\,&
\sum_{i=1}^{N}\sum_{k=0}^{\infty}\sum_{l}
\bigg\Vert 
\frac{1}{(\sigma^{k}_{i}+2\rho_{i})^{1/2}}
        \Big[ {1}_{m_{l}\times m_{l}}
                -\frac{\rho_{i}(\pcG_{i,l})^T\pcG_{i,l}}{\sigma^{k}_{i}+2\rho_{i}}
                \Big]
\notag\\[-3pt]&\hskip8mm
   \Big[
         2\rho_{i}(X^{k}_{i,l}-Z^{k}_{l})
        +\rho_{i}(\svZn^{k}_{i,l}-\svZp^{k}_{i,l})
        +\rho_{i}(\pcG_{i,l})^T(\pcG_{i}(X^{k+1,k}_{i,l})+Y^{k}_{i})
        +(f_{l})^T
        -\dvZp^{k}_{i,l}+\dvZn^{k}_{i,l}
        +(\pcG_{i,l})^T\dvG^{k}_{i}
        \Big]
\bigg\Vert^2.
\label{LPGSdeltaXkiSummed}
\end{align}
Adopt $\rho_{i}\geq 1$ as a reference here.
During the initial phase, the $\rho_{i}$-related terms,
$2\rho_{i}(X^{k}_{i,l}-Z^{k}_{l})
            +\rho_{i}(\svZn^{k}_{i,l}-\svZp^{k}_{i,l})
            +\rho_{i}(\pcG_{i,l})^T(\pcG_{i}(X^{k+1,k}_{i,l})+Y^{k}_{i})$
can differ significantly from zero.
To restrict the magnitude of the sum
on the left-hand side of \eqref{LPGSdeltaXkiSummed}, we take 
$(\sigma^{k}_{i}+2\rho_{i})^{1/2}\gg 2\rho_{i}$
to obtain $\sigma^{k}_{i} \gg (2\rho_{i})^2$.
Moreover, on the basis of \eqref{LPGSDualUpdateRules},
 \eqref{LPGSInitializationDual}, and
\eqref{LPGSProximalXZYAFOCSum2nd} through \eqref{LPGSdeltaXkiSummed}
and to have the sums in \eqref{LPGSProximalXZYAFOCSum2nd} finite, we propose
\begin{align}
\{\sigma^{k}_{i}, \tau^{k}\}  \gg (2\rho_{i})^2,\ \ 
\{\svZpcoef^{k}_{i}, \svZncoef^{k}_{i}, \gamma^{k}_{i}\}  \gg (\rho_{i})^2.
\label{LPGSControlParametersEstimated}
\end{align}
That is, the numerical values of the proximal control parameters are of
 orders of magnitude greater than $(\rho_{i})^2$.
A few relevant issues need to be studied further:
(a) 
The numerical values of $\rho_{i}$ should be relatively low to ensure 
that the values of the proximal control parameters are moderate.
(b)
Also, since $\rho_{i}$ affect $L^{0}$ directly, 
 pick their values adequately to have \eqref{LPGSProximalXZYAFOCSumInfty} hold.
(c)
Whether the relations in \eqref{LPGSControlParametersEstimated} can be relaxed
at higher $k$'s, 
considering that
$\{\Vert X^{k}_{i,l}-Z^{k}_{l}\Vert^2$,
            $\Vert \svZn^{k}_{i,l}-\svZp^{k}_{i,l}\Vert^2$,
            $\Vert \pcG_{i}(X^{k+1,k}_{i,l})+Y^{k}_{i}\Vert^2$,
            $\Vert \svZp^{k}_{i,l}\Vert^2$,
            $\Vert \svZn^{k}_{i,l}\Vert^2,\ldots\}$
are expected to be low at higher $k$'s.

To generate the bounded dual sequences
$\{\dvZp^{k}_{i}$, $\dvZn^{k}_{i}$, $\dvG^{k}_{i}\}_{k}$,
we apply \eqref{LPGSDualUpdateRules},
\eqref{LPGSdeltasvZnkSoluIterated}, \eqref{LPGSdeltasvZpkSoluIterated},
\eqref{LPGSdeltasvGkSoluIterated}, and the original primal update rules
for the slack variables to get
\begin{align}
&
\dvZp^{k+1}_{i}-\dvZp^{k}_{i}
=-\alpha^{k+1}_{\dvZp_{i}}\big(Z^{k+1}_{}-X^{k+1}_{i}+\svZp^{k+1}_{i}\big),\ \ 
\svZp^{k+1}_{i}-\svZp^{k}_{i}
=
-\frac{\rho_{i}}{\svZpcoef^{k}_{i}}
      \big(Z^{k+1}_{}-X^{k+1}_{i}+\svZp^{k+1}_{i}+\dvZp^{k}_{i}/\rho_{i}\big)
\notag\\[-3pt]&\hskip64mm\text{or}\
\svZp^{k+1}_{i}=
 P_{[0,u_{\dvZp_{i}}]}\Big(\svZp^{k}_{i}
   -\frac{\rho_{i}}{\svZpcoef^{k}_{i}}\big(Z^{k+1}_{}-X^{k+1}_{i}+\svZp^{k+1}_{i}
                        +\dvZp^{k}_{i}/\rho_{i}\big)\Big);   
\label{LPGSDeltasvdvZpEstimate}
\end{align} 
\vspace{-3mm}
\begin{align}
&
\dvZn^{k+1}_{i}-\dvZn^{k}_{i}
=-\alpha^{k+1}_{\dvZn_{i}}\big(X^{k+1}_{i}-Z^{k+1}_{}+\svZn^{k+1}_{i}\big),\ \ 
\svZn^{k+1}_{i}-\svZn^{k}_{i}
=
-\frac{\rho_{i}}{\svZncoef^{k}_{i}}
      \big(X^{k+1}_{i}-Z^{k+1}_{}+\svZn^{k+1}_{i}+\dvZn^{k}_{i}/\rho_{i}\big)
\notag\\[-3pt]&\hskip64.5mm\text{or}\
\svZn^{k+1}_{i}=
 P_{[0,u_{\dvZn_{i}}]}\Big(\svZn^{k}_{i}
   -\frac{\rho_{i}}{\svZncoef^{k}_{i}}\big(X^{k+1}_{i}-Z^{k+1}_{}+\svZn^{k+1}_{i}
                        +\dvZn^{k}_{i}/\rho_{i}\big)\Big);   
\label{LPGSDeltasvdvZnEstimate}
\end{align} 
and
\begin{align}
&
\dvG^{k+1}_{i}-\dvG^{k}_{i}
=-\alpha^{k+1}_{\dvG_{i}}\big(\pcG_{i}(X^{k+1}_{i})+Y^{k+1}_{i}\big),\ \
Y^{k+1}_{i}-Y^{k}_{i}
=
-\frac{\rho_{i}}{\gamma^{k}_{i}}
      \big(\pcG_{i}(X^{k+1}_{i})+Y^{k+1}_{i}+\dvG^{k}_{i}/\rho_{i}\big)
\notag\\[-3pt]&\hskip70mm\text{or}\ 
Y^{k+1}_{i}=
 P_{[0,u_{Y_{i}}]}\Big(Y^{k}_{i}
   -\frac{\rho_{i}}{\gamma^{k}_{i}}\big(\pcG_{i}(X^{k+1}_{i})+Y^{k+1}_{i}+\dvG^{k}_{i}/\rho_{i}\big)\Big).     
\label{LPGSDeltasvdvGiEstimate}
\end{align} 
The above three sets of paired relations are possible 
owing to the penalty ($\rho_{i}>0$) employed.
To proceed further, we partition $\{\pcG_{i}(Z)\leq 0\}$ into two subgroups:
one is $\{\pm \pcH_{i}(Z)\leq 0\}$ from the original problem~\ref{LPGSPrimalProblemConstraintSeto}, 
denoted as $\{\pcH_{i}\}$ below;
the other is $\{\pcG_{i}(Z)\leq 0\}$ proper
from the original problem~\ref{LPGSPrimalProblemConstraintSeto}, 
denoted as $\{\pcG_{i}/\pcH\}$.

Motivated by the need to produce the bounded dual sequences 
and the mathematical structures of the paired expressions in
\eqref{LPGSDeltasvdvZpEstimate}, \eqref{LPGSDeltasvdvZnEstimate},
and \eqref{LPGSDeltasvdvGiEstimate}, we propose
\begin{align}
 \alpha^{k+1}_{\dvZp_{i}} \leq \frac{\rho_{i}}{\svZpcoef^{k}_{i}},\ \
 \alpha^{k+1}_{\dvZn_{i}} \leq \frac{\rho_{i}}{\svZncoef^{k}_{i}},\ \
 \alpha^{k+1}_{\dvG_{i}}  \leq \frac{\rho_{i}}{\gamma^{k}_{i}}\text{ (for $\{\pcH_{i}\}$)},\ \
 \alpha^{k+1}_{\dvG_{i}}  \geq \frac{\rho_{i}}{\gamma^{k}_{i}}\text{ (for $\{\pcG_{i}/\pcH\}$)}.
\label{LPGSDualCoeffsVsControlParas}
\end{align}
Along with the other initial and parameter values estimated above,
an appropriate combination of 
$\{\alpha^{k+1}_{\dvZp_{i}}$, $\alpha^{k+1}_{\dvZn_{i}}$, $\alpha^{k+1}_{\dvG_{i}}$,
${\rho_{i}}/{\svZpcoef^{k}_{i}}$, 
${\rho_{i}}/{\svZncoef^{k}_{i}}$, 
${\rho_{i}}/{\gamma^{k}_{i}}\}$ obeying \eqref{LPGSDualCoeffsVsControlParas}
and $\{\dvZp^{0}_{i}$,  $\dvZn^{0}_{i}$, $\dvG^{0}_{i}$,
$\svZp^{0}_{i}$,  $\svZn^{0}_{i}$, $Y^{0}_{i}\}$ set above
is expected to generate the dual sequences bounded. 
To examine the contents and consequences of \eqref{LPGSDeltasvdvZpEstimate}
through \eqref{LPGSDualCoeffsVsControlParas} in some detail, let us recall  
the limits from Lemma~\ref{LPGSPDDeltaLimitsp}(iii)(iv),
\begin{align}
\dvZp^{\infty}_{i}\geq 0,\ \
\svZp^{\infty}_{i}= 0;\ \
\dvZn^{\infty}_{i}\geq 0,\ \
\svZn^{\infty}_{i}= 0;\ \
\dvG^{\infty}_{i}\geq 0,\ \
 Y^{\infty}_{i}=-\pcG_{i}(X^{\infty}_{i})\geq 0,\ \
\big\langle \dvG^{\infty}_{i}, Y^{\infty}_{i}\big\rangle=0.
\label{LPGSSequLimits}
\end{align} 

Some rationals for the first relation of 
\eqref{LPGSDualCoeffsVsControlParas} are as follows.
(a) 
$\dvZp^{\infty}_{i}\geq 0$ and $\svZp^{\infty}_{i}=0$
from \eqref{LPGSSequLimits}; they suggest possibly that the rate 
of $\{\dvZp^{k}_{i}\}_{k}$ should be lower than that of $\{\svZp^{k}_{i}\}_{k}$
under \eqref{LPGSInitializationDual}.
(b)
According to \eqref{LPGSDeltasvdvZpEstimate},
$\alpha^{k+1}_{\dvZp_{i}}$ and $\rho_{i}/\svZpcoef^{k}_{i}$
control directly the rates of $\{\dvZp^{k}_{i}\}_{k}$
and $\{\svZp^{k}_{i}\}_{k}$, respectively.
(c)
Combination of (a) and (b)
results in $\alpha^{k+1}_{\dvZp_{i}} \leq \rho_{i}/\svZpcoef^{k}_{i}$.
(d) 
Together with other parameter and initial values,
an appropriate combination of low $\alpha^{k+1}_{\dvZp_{i}}$ and high $\dvZp^{0}_{i}$
may help generate the desired positive sequence
$\{\dvZp^{k}_{i}\}_{k}$ bounded through (\ref{LPGSDeltasvdvZpEstimate}$)_1$.
(e) 
On the one hand,
relatively higher $\{\rho_{i}/\svZpcoef^{k}_{i}\}_{k}$
makes $\{\svZp^{k}_{i}\}_{k}$ decreasing faster
        according to (\ref{LPGSDeltasvdvZpEstimate}$)_{2}$,
the resulting relatively lower $\{\svZp^{k}_{i}\}_{k}$
in turn makes $\{\dvZp^{k}_{i}\}_{k}$ decreasing slower,
        according to (\ref{LPGSDeltasvdvZpEstimate}$)_{1}$.
On the other hand,
relatively lower $\{\alpha^{k}_{\dvZp_{i}}\}_{k}$
slows down the decrease of $\{\dvZp^{k}_{i}\}_{k}$ 
(or even increasing if $Z^{k+1}_{}-X^{k+1}_{i}+\svZp^{k+1}_{i}<0$ or the like)
        according to (\ref{LPGSDeltasvdvZpEstimate}$)_{1}$,
 the resultant relatively higher $\{\dvZp^{k}_{i}\}_{k}$
may tend to make $\{\svZp^{k}_{i}\}_{k}$ decreasing faster toward zero,
according to (\ref{LPGSDeltasvdvZpEstimate}$)_{2,3}$,
especially through the term of $-(\rho_{i}/\svZpcoef^{k}_{i})(\dvZp^{k}_{i}/\rho_{i})$.
Thus, there is a favorable feedback mechanism to yield the desired sequences,
$\{\dvZp^{k}_{i},\svZp^{k}_{i}\}_{k}$.

Justifications parallel to the above can be offered for 
the second relation of \eqref{LPGSDualCoeffsVsControlParas} too.

Regarding the third relation for the constraints of $\{\pcH_{i}\}$
in \eqref{LPGSDualCoeffsVsControlParas},
the limits~\eqref{LPGSSequLimits} give
$\{\dvG^{\infty}_{i,j}\geq 0$, $Y^{\infty}_{i,j}=0\}$ or
$\{\dvG^{\infty}_{i,j}> 0$, $Y^{\infty}_{i,j}=0\}$ most possibly,
as required by the stationarity conditions of Lemma~\ref{LPGSPDDeltaLimitsp}(iv).
Reasons similar to (b) through (e) above 
can be used here to argue for $\alpha^{k+1}_{\dvG_{i}} \leq \rho_{i}/\gamma^{k}_{i}$.

We now discuss the fourth relation for $\{\pcG_{i}/\pcH\}$
 in \eqref{LPGSDualCoeffsVsControlParas}.
(a) 
The limits~\eqref{LPGSSequLimits} give
$\{\dvG^{\infty}_{i}\geq 0$, $Y^{\infty}_{i}=-\pcG_{i}(X^{\infty}_{i})\geq 0$,
$\big\langle \dvG^{\infty}_{i}, Y^{\infty}_{i}\big\rangle=0\}$; 
if $Y^{\infty}_{i,j}> 0$,
then $\dvG^{\infty}_{i,j}=0$;
 this expectedly occurs to the majority of $\{\pcG_{i}(X^{\infty}_{i})\leq 0\}$ inactive
in the original problem~\ref{LPGSPrimalProblemConstraintSeto},
that is, the limits of $\{\dvG^{\infty}_{i,j}=0$, $Y^{\infty}_{i,j}>0\}$
need to be produced most possibly.
(b)
To produce the limits,
the structures of equations~\eqref{LPGSDeltasvdvGiEstimate} suggest
 $\alpha^{k+1}_{\dvG_{i}}\geq\rho_{i}/\gamma^{k}_{i}$ with large $\gamma^{k}_{i}$ taken,
 while $\{\dvG^{0}_{i}$, $Y^{0}_{i}\}$ assigned according to 
 \eqref{LPGSInitializationDual} and \eqref{LPGSDualControlParametersEstimated}.
(c)
 This strict inequality tends to have $Y^{k+1}_{i}\geq 0$
 through the term of $\dvG^{k}_{i}/\rho_{i}$
 in (\ref{LPGSDeltasvdvGiEstimate}$)_{2}$ as explained below.  
  Greater $\alpha^{k+1}_{\dvG_{i}}$ tends to make $\{\dvG^{k}_{i}\}_{k}$ lower 
by (\ref{LPGSDeltasvdvGiEstimate}$)_{1}$,
which in turn tends to make $\{Y^{k}_{i}\}_{k}$ higher 
by $\dvG^{k}_{i}/\rho_{i}$ via (\ref{LPGSDeltasvdvGiEstimate}$)_{2}$, 
yielding the desired nonnegative and high sequence $\{Y^{k}_{i}\}_{k}$ 
simply via (\ref{LPGSDeltasvdvGiEstimate}$)_{2}$,
instead of (\ref{LPGSDeltasvdvGiEstimate}$)_{3}$.
(d) 
Next, suppose that together with $Y^{0}_{i}$,
        $\gamma^{k}_{i}$ is sufficiently large to generate
a nonnegative and high sequence $\{Y^{k}_{i}\}_{k}$ through the explicit expression
of $Y^{k+1}_{i}$, (\ref{LPGSDeltasvdvGiEstimate}$)_{2}$.
If the relatively large $\alpha^{k+1}_{\dvG_{i}}$ and high $Y^{k+1}_{i}$
yield $\dvG^{K+1}_{i,j}<0$ from 
(\ref{LPGSDeltasvdvGiEstimate}$)_{1}$/(\ref{LPGSDualUpdateRules}$)_{3}$ 
for some $K\in\mathbb{N}$,
there is the necessity to bound $\{\dvG^{k}_{i}\}_{k}$ explicitly from below 
by zero,
and the update model (\ref{LPGSDualUpdateRules}$)_{3}$ for $\{G_{i}/H\}$
is modified as follows.
Instead of a scalar,
the dual coefficient $\alpha^{k+1}_{\dvG_{i}}$ is extended to be a
semi-positively definite diagonal matrix defined through
\vspace{-2mm}
\begin{align}
 \alpha^{k+1}_{\dvG_{i},j}=
\l\{\begin{array}{ll}
\alpha_{\dvG_{i}}, & 
        \text{if $\tilde{\dvG}^{k+1}_{i,j}:=\dvG^{k}_{i,j} -\alpha_{\dvG_{i}} \big(\pcG_{i}(X^{k+1}_{i})+Y^{k+1}_{i}\big)_{j}
                        \in [0, u_{\dvG_{i},j}]$,}\\[3pt]
0, & \text{else},
    \end{array}\r. 
\label{LPGSDualUpdateRulesdvWBoundsj}
\end{align}
where $\alpha_{\dvG_{i}}$ is a positive scalar constant.
Accordingly,
the relation $\sum_{k=1}^{\infty}\Vert \pcG_{i}(X^{k}_{i})+Y^{k}_{i} \Vert^2<\infty$
in Lemma~\ref{LPGSPDSumBoundedp}(i) is modified to
$\sum_{k=0}^{\infty}\big\langle\alpha^{k+1}_{\dvG_{i}}(\pcG_{i}(X^{k+1}_{i})+Y^{k+1}_{i}),
                        \pcG_{i}(X^{k+1}_{i})+Y^{k+1}_{i}\big\rangle<\infty$,
which does not imply a definite consequence about $\{\pcG_{i}(X^{k}_{i})+Y^{k}_{i}\}_{k}$
 for $\{G_{i}/H\}$. 
Furthermore, the inequality $\alpha^{k+1}_{\dvG_{i}}\geq\rho_{i}/\gamma^{k}_{i}\text{ (for $\{\pcG_{i}/\pcH\}$)}$
in \eqref{LPGSDualCoeffsVsControlParas} needs to be modified to
 $\alpha^{}_{\dvG_{i}}\geq\rho_{i}/\gamma^{k}_{i}\text{ (for $\{\pcG_{i}/\pcH\}$)}$.
We now need to resolve the issues whether 
$\lim_{k\rightarrow\infty}\dvG^{k}_{i}$ exists,
$\lim_{k\rightarrow\infty}(\pcG_{i}(X^{k}_{i})+Y^{k}_{i})=0$, and 
$\lim_{k\rightarrow\infty}\langle \dvG^{k}_{i},Y^{k}_{i}\rangle=0$ 
for $\{G_{i}/H\}$ under \eqref{LPGSDualUpdateRulesdvWBoundsj}.    
(e)
To this end, we consider two bounded inequalities,
$\sum_{k=0}^{\infty}\big\langle\alpha^{k+1}_{\dvG_{i}}(\pcG_{i}(X^{k+1}_{i})+Y^{k+1}_{i}),
                        \pcG_{i}(X^{k+1}_{i})+Y^{k+1}_{i}\big\rangle<\infty$
and $\sum_{k=0}^{\infty}\Vert Y^{k+1}_{i}-Y^{k}_{i} \Vert^2
=\sum_{k=0}^{\infty}\Vert ({\rho_{i}}/{\gamma^{k}_{i}})
                (\pcG_{i}(X^{k+1}_{i})+Y^{k+1}_{i}+\dvG^{k}_{i}/\rho_{i})\Vert^2<\infty$,
the latter holds under Lemma~\ref{LPGSPDSumBoundedp}(ii) and
(\ref{LPGSDeltasvdvGiEstimate}$)_2$ argued above.
Supposing that $\{\gamma^{k}_{i}\}_{k}$ is bounded from above by a finite value
and applying (\ref{LPGSDualUpdateRules}$)_{3}$ and \eqref{LPGSDualUpdateRulesdvWBoundsj} to the two inequalities,
we obtain
\begin{align}
&
 \sum_{\alpha^{k+1}_{\dvG_{i},j}=\alpha_{\dvG_{i}}}\ss
        \big\vert \big(\pcG_{i}(X^{k+1}_{i})+Y^{k+1}_{i}\big)_{j}\big\vert^2<\infty,\ 
 \sum_{\alpha^{k+1}_{\dvG_{i},j}=\alpha_{\dvG_{i}}}\ss
        \big\vert \dvG^{k}_{i,j} \big\vert^2<\infty,\
\sum_{\alpha^{k+1}_{\dvG_{i},j}=0}\ss
        \big\vert \big(\pcG_{i}(X^{k+1}_{i})+Y^{k+1}_{i}+\dvG^{k}_{i}/\rho_{i}\big)_{j} \big\vert^2
                <\infty,\ \forall j.
\label{LPGSsdvGConditionsj}
\end{align}
Consider the $j$-th components,
three cases need to be discussed as follows.
\begin{enumerate}
\item[i.]
$\exists K\in\mathbb{N}$ such that $\{\alpha^{k+1}_{\dvG_{i},j}=\alpha_{\dvG_{i}}\}_{k\geq K}$.
\eqref{LPGSsdvGConditionsj} reduce to
\begin{align*}
&
 \sum_{k=K}^{\infty}
        \big\vert \big(\pcG_{i}(X^{k+1}_{i})+Y^{k+1}_{i}\big)_{j}\big\vert^2<\infty,\ \ 
 \sum_{k=K}^{\infty}
        \big\vert \dvG^{k}_{i,j} \big\vert^2<\infty,
\end{align*}
which, together with Lemma~\ref{LPGSPDAccumulationPointp}, gives
$\dvG^{\infty}_{i,j}= 0$, 
$Y^{\infty}_{i,j}\geq 0$,
$(\pcG_{i}(X^{\infty}_{i})+Y^{\infty}_{i})_{j}=0$, 
$\big\langle \dvG^{\infty}_{i,j}, Y^{\infty}_{i,j}\big\rangle=0$.
\item 
[ii.]
$\exists K\in\mathbb{N}$ such that $\{\alpha^{k+1}_{\dvG_{i}}=0\}_{k\geq K}$.
(\ref{LPGSsdvGConditionsj}$)_{3}$,
\eqref{LPGSDualUpdateRulesdvWBoundsj},
 (\ref{LPGSDualUpdateRules}$)_{3}$,
and $\{\dvG^{k}_{i,j}\}_{k}$ bounded from below by zero give
\begin{align*}
&
\sum_{k=K}^{\infty}
        \big\vert \big(\pcG_{i}(X^{k+1}_{i})+Y^{k+1}_{i}+\dvG^{k}_{i}/\rho_{i}\big)_{j} \big\vert^2
                <\infty,\ \ 
\dvG^{k+1}_{i,j}=\dvG^{k}_{i,j}\geq 0,\ \
\dvG^{k}_{i,j}-\alpha_{\dvG_{i}} \big(\pcG_{i}(X^{k+1}_{i})+Y^{k+1}_{i}\big)_{j}<0,
 \ \forall k\geq K.
\end{align*}
Together with Lemma~\ref{LPGSPDAccumulationPointp}, they result in
$\dvG^{\infty}_{i,j}= 0$, 
$Y^{\infty}_{i,j}\geq 0$,
$(\pcG_{i}(X^{\infty}_{i})+Y^{\infty}_{i})_{j}=0$, 
$\big\langle \dvG^{\infty}_{i,j}, Y^{\infty}_{i,j}\big\rangle=0$.
\item[iii.]
Consider the mixed cases between Case i and Case ii.
It follows from (\ref{LPGSsdvGConditionsj}$)_2$ 
and $\dvG^{k+1}_{i,j}=\dvG^{k}_{i,j}$ under $\alpha^{k+1}_{\dvG_{i},j}=0$
that $\lim_{k\rightarrow\infty}\dvG^{k}_{i,j}=0$.
Next, (\ref{LPGSsdvGConditionsj}$)_{1,3}$ yield
 $\lim_{k\rightarrow\infty}\big(\pcG_{i}(X^{k}_{i})+Y^{k}_{i}\big)_{j}=0$.
\end{enumerate}
It then follows from the three cases that
$\dvG^{\infty}_{i}= 0$, 
$Y^{\infty}_{i}\geq 0$,
$\pcG_{i}(X^{\infty}_{i})+Y^{\infty}_{i}=0$, and
$\big\langle \dvG^{\infty}_{i}, Y^{\infty}_{i}\big\rangle=0$
for $\{\pcG_{i}/\pcH\}$.
Thus, the fourth relation
 in \eqref{LPGSDualCoeffsVsControlParas} modified, 
 $\alpha^{}_{\dvG_{i}}\geq\rho_{i}/\gamma^{k}_{i}\text{ (for $\{\pcG_{i}/\pcH\}$)}$
 is justified.
 We notice that $\dvG^{\infty}_{i}= 0$ for $\{\pcG_{i}/\pcH\}$
 is allowable, 
considering the roles played by $\dvZp^{\infty}_{i}$ and $\dvZn^{\infty}_{i}$
in the stationarity conditions of Lemma~\ref{LPGSPDDeltaLimitsp}(iv).

\section{Discussion and Summary} 
 \label{sec:Summary}
 
 The present study develops an algorithm composed of 
 the update rules~\eqref{LPGSProximalXilA},
\eqref{LPGSProximalZl},  
\eqref{LPGSProximalsvZpilA}, \eqref{LPGSProximalsvZnilA},
\eqref{LPGSProximalYilA},
and \eqref{LPGSDualUpdateRules}/\eqref{LPGSDualUpdateRulesdvWBoundsj}
 for distributed computing of linear programming
 of huge-scales. Global consensus with single common variable, multiblocks, and
augmented Lagrangian are adopted.
The global consensus is used to partition the constraints of equality and inequality into 
$N$ consensus blocks,
 and the subblocks of each consensus block are used to partition the primal variables into $M$ sets
 of disjoint subvectors. 
 The global consensus constraints of equality and other constraints are converted into
 the extended constraints of equality involving slack variables
 in order to help the feasibility and initialization of the algorithm.
 The block-coordinate Gauss-Seidel method, the proximal point method,
 and ADMM are used to update the primal variables, 
and descent models are used to update the dual.
 Convergence of the algorithm to optimal solutions
 is argued and the rate of convergence, $O(1/k^{1/2})$ is estimated, 
 under feasibility of the algorithm and boundedness of the dual sequences supposed.
 The issues of how to make the algorithm feasible and the dual sequences bounded are discussed.

 A few issues listed below need to be clarified and/or explored further.
 
First, 
the original global consensus constraints of equality, $\{X_{i}-Z=0\}$
and others, $\{\pcH(Z)=0$, $\pcG(Z)\leq 0\}$ are replaced equivalently
by the extended constraints of equality involving slack variables;
the flexibility of setting initial values for the slack variables helps
make the algorithm feasible and dual sequences bounded,
though the employment of the slack variables increases 
computational size of the resultant optimization problem.
   
Second,
the dual descent update models \eqref{LPGSDualUpdateRules}
and the extended \eqref{LPGSDualUpdateRulesdvWBoundsj} are adopted.
Their introduction is motivated by the mathematical structure of 
Lemma~\ref{LPGSProximalFOCp} and the satisfaction of the constraints
of the global consensus problem \eqref{LPGSPrimalProblemGlobalConsensusRing} 
subject to \eqref{LPGSConstraintsithCBSlacked}
by the limits of the primal sequences,
their merit is supported by the results stated in 
Theorems~\ref{LPGSConvergenceT} and \ref{LPGSConvgRate}.
We observe that the dual descent updates are a consequence 
of the primal updates leading to Lemma~\ref{LPGSProximalFOCp},
and we may view the development as algorithmic modeling in order to determine solutions
to the primal problem~\eqref{LPGSPrimalProblemOriginal},  
similar to constitutive modeling in continuum mechanics.
It is interesting to notice the link between the present and that
 discussed in \cite{SunSun2024} which concerns  dual descent augmented Lagrangian method
 involving highly nonconvex constraints,
 while the present deals with linear programming.
 
Third, 
the penalty terms in the augmented Lagrangian avail the paired expressions
of dual and slack variables
presented in \eqref{LPGSDeltasvdvZpEstimate}, \eqref{LPGSDeltasvdvZnEstimate},
and \eqref{LPGSDeltasvdvGiEstimate}.
These paired expressions suggest how to generate bounded dual sequences
and the desired limits through adequate choices of the control parameters and initial values.
  
Fourth,
as analyzed in Sec.~\ref{sec:Initialization},
feasibility of the algorithm and boundedness of the dual sequences
 need to be ensured
through 
the initial values of the primal and dual sequences,
the control parameter values of
$\{\rho_{i}$, $\sigma^{k}_{i}$, $\tau^{k}$, 
$\svZpcoef^{k}_{i}$, $\svZncoef^{k}_{i}$, $\gamma^{k}_{i}$, 
$\alpha^{k+1}_{\dvZp_{i}}$, $\alpha^{k+1}_{\dvZn_{i}}$, $\alpha^{k+1}_{\dvG_{i}}\}$,
 and preconditioning of the constraints of $\{\pcG_{}(Z)\leq 0$, $\pcH_{}(Z)= 0\}$ 
 via scaling-down. 
 These issues are to be studied further to specify the values
 and optimize the computation. 
  
Fifth,
the rate of convergence of the algorithm, $O(1/k^{1/2})$ is estimated crudely and
a better treatment should be explored.

Sixth,
we have discussed the behavior of the dual model with built-in bound for $\{\pcG_{i}|\pcH\}$,
\eqref{LPGSDualUpdateRulesdvWBoundsj},
supposing that (\ref{LPGSDeltasvdvGiEstimate}$)_2$ holds for all $k$
(and its resultant $\{y^{k}_{i}\}_{k}$ is nonnegative).
Have yet to analyze the model's behavior to deal with $Y^{k+1}_{i,j}=0$
and $\partial_{Y_{i,j}}\Indicator_{[0, u_{Y_{i},j}]}(Y^{k+1}_{i,j})$ at finite $k$.
Have yet to discuss the possibility of dual models with built-in bound for 
$\{\dvZp^{k}_{i}$, $\dvZn^{k}_{i}$, $\dvG^{k}_{i}\ \text{(for $\{\pcH_{i}$\})}\}$.
 
Finally,
 how to extend 
the present formulation beyond the linear programming represented 
by \eqref{LPGSPrimalProblemOriginal} and \eqref{LPGSPrimalProblemConstraintSeto}?
One example is to resolve the issue of non-convergence (of direct extension of ADMM
for multi-block convex optimizations) raised in \cite{Chenetal2016}
 for aggregatively computable and convex objective functions, with the aid of 
the present treatment.
The other is to extend the primal problem to the mathematical form of
$\{f(S), H(S,T)=0, H'(T,Z)=0, G(S,T,Z)\leq 0 \}$.
Here, the dimension of $S$ is much smaller than that of $T$
and the latter smaller than that of $Z$.
The third example is how to extend the treatment to aggregatively computable
and quadratic constraints of inequality, convex and difference of convex.
The ultimate goal is to develop an algorithm for CFD-optimization coupled problems
of inhomogeneous turbulent flows.

\bibliography{DCHSLP.bib}

\end{document}